\documentclass[letterpaper]{article}
\usepackage[a4paper,top=0.75in,bottom=1.05in,left=0.75in,right=0.75in,marginparwidth=0.75cm]{geometry}
\usepackage{amsmath}
\usepackage{amssymb}
\usepackage{mathtools} 
\usepackage{enumerate}
\usepackage{graphicx}
\usepackage{float}
\usepackage[superscript, biblabel, nomove]{cite}
\usepackage{color}
\usepackage{url}
\usepackage{subcaption}
\usepackage{caption}
\usepackage{caption}
\usepackage{hyperref}
\hypersetup{pdfstartview={XYZ null null 1.00}}
\usepackage{amsmath}
\usepackage{amssymb}
\usepackage{mathtools} 
\usepackage{enumerate}
\usepackage{graphicx}
\usepackage{float}

\usepackage{color}
\usepackage{url}
\usepackage{subcaption}
\usepackage{caption}
\usepackage{hyperref}
\usepackage{amssymb}
\usepackage{amsmath}
\usepackage{amsthm}
\usepackage{float}
\usepackage{siunitx}
\usepackage{booktabs}
\usepackage[english]{babel}
\usepackage[nottoc]{tocbibind}
\usepackage{lipsum}
\usepackage{booktabs}
\usepackage{siunitx}
\usepackage{makecell}
\usepackage{changepage}
\usepackage{orcidlink}
\usepackage{hyperref}
\usepackage{authblk}
\usepackage{blindtext}
\usepackage{datetime}
\definecolor{DarkBlue}{rgb}{0.00,0.00,0.50}
\hypersetup{colorlinks=true,
            linkcolor=blue,
            filecolor=magenta,
            urlcolor=cyan,
            citecolor=DarkBlue,
            pdftitle={MSc UoR Computer Science Report},
            }
\begin{document}
\title{An Efficient Approach to Fractional Analysis for Non-Linear Coupled Thermo-Elastic Systems}
\date{}
\author
{Qasim Khan }
\affil{{ Department of Mathematics and Information Technology,\\  The Education University of Hong Kong,\\ 10 Lo Ping Road, Tai Po, N.T,\\ Hong Kong} }
\affil{ {qasimkhan@s.eduhk.hk }}
\maketitle	
\abstract{Nonlinear thermoelastic systems play a crucial role in understanding thermal conductivity, stresses, elasticity, and temperature interactions. This research focuses on finding solutions to these systems in their fractional forms, which is a significant aspect of the study. We consider various proposed models related to fractional thermoelasticity and derive results through sophisticated methodologies. Numerical simulations are conducted for both fractional and integer order thermoelastic coupled systems, with results presented in tables and graphs. The graphs indicate a close correspondence between the approximate and exact solutions. The solutions obtained demonstrate convergence for both fractional and integer order problems, ensuring accurate modeling. Furthermore, the tables confirm that greater accuracy can be achieved by increasing the number of terms in the series of solutions. }

\section{Introduction}\label{s1}
Several scientific fields, such as fluid mechanics, solid-state physics, and plasma physics, deal with the  non-linear one-dimensional coupled systems of the  thermo-elasticity. Thermo-electricity problems have received much attention due to their various applications in different branches. Thermoelasticity is an area in applied mechanics that studies strains and stresses as well as heat transfer theory. Deformation of solids can happen as a result of heat flow, and temperature changes can happen as a result of strain fields\cite{thermothermal1}. Interaction between mechanical and thermal domains can be studied using coupled linear and non-linear thermo-elasticity  systems. The study of these principles have recently aroused the interest of a number of mathematicians working in a variety of research topics. Famous scientists, mathematicians, and engineers have been influenced by the proof of irrational physical behavior, demonstrated by elastic deformations caused by temperature stresses. Duhame presented the classical uncoupled thermoelasticity theory in \cite{thermothermal2}.
 Maurice A. Biot extended this theory to coupled thermoelasticity in \cite{thermothermal3}.

In this connection, many scientists examine scientific phenomena' behaviors and dynamical features, using the basic concepts and theory of fractional calculus (FC)\cite{thermo5}. Recently, [Iqbal 2023 \cite{thermoNaveed}] used a transformation-based homotopy perturbation method with the Caputo operator to illustrate fractional order singular and non-singular thermoelastic coupled systems. In the same regard, Rab et al. used a Laplace residual power series method [LRPSM] for the analytical solution of the thermoelastic coupled systems ref. \cite{thermoShabaz}. Similarly, Hassan Eltayeb Gadain combined the double Laplace transform with the decomposition method to solve integer order nonlinear singular one-dimensional thermo-elasticity coupled system ref.\cite{thermo32,thermo32a}.

In recent decades, fractional partial differential equations (FPDEs) have been considered one of the most powerful tools for modeling and have attracted many researchers to the topic of fractional calculus (FC). Because most of the non-linear fractional-order problems do not have exact solutions, and thus analytical and numerical methodologies have been employed to examine their approximation solutions.  In this regard, the following efficient and novel techniques such as, ADM \cite{thermo17}, VIM \cite{thermo16}, q-HAM \cite{thermo18}, operational matrix method \cite{thermo20}, the fast element-free Galerkin method \cite{thermo19}, new analytical approximate method (NAAM) \cite{thermo21}, iterative Laplace  transform method \cite{thermosacacaedc22}, New iterative Aboodh transform \cite{thermo22}, operational calculus scheme\cite{thermo24},  Homotopy transform technique (NIATM)\cite{thermo23}, Elzaki transform decomposition method (ETDM) \cite{thermo25}, iterative reproducing kernel method \cite{thermo27} and multistage differential transformation method \cite{thermo26} have recently been developed and successfully applied to different fractional partial differential equation.

 Integral transformations are the most influential mathematical approaches for solving differential equations, partial differential equations, integro-differential equations, partial integro-differential equations, delay differential equations, population growth problems and fractional order related problems. The well-known transformations are the Laplace, Fourier, Elzaki, and Sumudu transforms. A novel integral transform known as Aboodh transform is similar to the Laplace transform and other integral transforms in the time domain ${{{t}}}>0$ such as the  Natural transform, Elzaki transform, and Sumudu transform, \cite{thermoAT,thermoAT1,thermoATF}. George Adomian, between 1970s and 1990s developed the Adomian decomposition method (ADM). ADM is an analytical method that is used for the solution of different dynamical systems. To control the non-linearity in a given system, different polynomials are introduced such as He, Orthogonal, and Jafari polynomials. For all types of non-linearity, Adomian \cite{thermo12X,thermoy,thermoz} formally presented formulas that can generate Adomian polynomials.

In the present article, the Aboodh transformation along with the decomposition method is applied to the solutions of the non-linear coupled systems of thermo-elasticity. The derivative is taken in the Caputo sense. The Aboodh transform is the generalization of the Laplace transform which converts the targeted problems into their simplest form. The Adomian decomposition method is then implemented to achieve the desired solution. The algorithm of the proposed method was developed using MAPLE software. The obtained results are shown by using graphs and tables for both fractional and integer order problems. It is observed that as the degree of the polynomial increases, the accuracy of the suggested method increases. The non-linearity of the problem is handled with the help of Adomian polynomials. The procedure to find the Adomian polynomials is tough but quite stable and accurate. The fractional-order solutions are calculated successfully and are convergent toward the integer-order solutions.

A short summary of the rest of the paper is provided below. In section \ref{s2} we discuss some fundamental preliminary concepts, in section \ref{s3}, the basic methodology is discussed for non-linear one-dimensional thermo-elasticity coupled system, and section \ref{s4} discusses some test models that demonstrate the effectiveness of ATDM,  and section \ref{s6} is the conclusion of the current manuscript.
\section{Basic preliminaries concepts }\label{s2}
In this section, we will discuss fundamental preliminaries that can be used in the rest of the units.
\subsection{Definition Michele Caputo, 1967,\cite{thermocaputo}}
$^{C}D_{{{{t}}}}^{{\beta}}$ is the Caputo derivative of a fractional order $\beta$ for the function ${{u}}({{{{x}}}},{{{t}}})$ and is defined as
\begin{equation}\label{p1}
  ^{C}D_{{{{t}}}}^{{\beta}}{{u}}({{{{x}}}},{{{t}}})=J_{{{{t}}}}^{\rho-{\beta}}{{u}}^{{\rho}}({{{{x}}}},{{{t}}}),  \ \ \rho-1<\beta\leq \rho,\ \ \rho\in N \ \ {{{t}}}>0.
\end{equation}
where   $J_{{{{t}}}}^{\gamma}$ is the Riemann-Liouville integral for fractional-order $\beta$  of the function ${{u}}({{{{x}}}},{{{t}}})$  and  define as
\begin{equation}\label{p11}
  J_{{{{t}}}}^{\beta}{{u}}({{{{x}}}},{{{t}}})=\frac{1}{\Gamma(\beta)}\int_{0}^{{{{t}}}}({{{t}}}-\kappa){{u}}({{{{x}}}},{{{t}}})d{{{t}}},\ \ \ \beta>0,\  {{{t}}}>0, \ \kappa\geq0,
\end{equation}
by considering that the above integral exists.
\subsection[Abdon Atangana]{Proposition} If $\rho-1<\beta\leq \rho$, $\rho \in \mathcal{N}^+ $ and $^{C}D_{{{{t}}}}^{{\beta}}{{u}}({{{{x}}}},{{{t}}})$ exists.\\
Then following properties hold in terms of Caputo sense:
 \begin{equation*}
 \begin{split}
      {\mathrm{\lim}_{\beta \rightarrow \rho}} ^{C}D_{{{{t}}}}^{{\beta}}{{u}}({{{{x}}}},{{{t}}})&=\frac{\partial^{\rho}}{\partial {{t}}^{\rho}} {{u}}({{{{x}}}},{{{t}}}),\\
 {\mathrm{\lim}_{\beta \rightarrow \rho-1}} ^{C}D_{{{{t}}}}^{{\beta}}{{u}}({{{{x}}}},{{{t}}})&=\frac{\partial^{\rho-1}}{\partial {{t}}^{\rho-1}} {{u}}({{{{x}}}},{{{t}}})-\frac{\partial^{\rho-1}}{\partial {{t}}^{\rho-1}} {{u}}({{{{x}}}},{{0}}).\\
\end{split}
 \end{equation*}
 \textbf{Proof:}
 Proof can be found in (Abdon Atangana, 2017  in Fractional Operators \cite{thermoAtangana2017}) chapter 5 section 3 on page 82.
\subsection{Definition [K Suliman Aboodh, 2013]}
Khalid Suliman Aboodh developed  a new integral transform known as the Aboodh transform for a function ${{u}}({{{x}}},{{{t}}})$ under a set $\mathbf{Z}$, and is defined as\cite{thermoAT,thermoAT1,thermoATF}
\begin{equation}
\mathbf{Z}=\{{{u}}({{{x}}},{{{t}}}):\exists M, \Bbbk_{1}, \Bbbk_{2}>0, |{{u}}({{{x}}},{{{t}}})|<M\},
\end{equation}
where constant M must be a finite number, $\Bbbk_1$,$\Bbbk_2$  may or mayn't be finite. \\
The Aboodh transform (AT) denoted by the operator $\mathbf{A}$ (.)  in this paper and is define by\\
\begin{equation}
\mathbf{A}\left[{{{u}}}({{{x}}},{{{t}}})\right]={k}({s})=\frac{1}{s}\int_{0}^{\infty}{{{u}}}({{{x}}},{{{t}}})e^{-s{{{t}}}}{d}{{{t}}},{{{t}}}\geq0,{\Bbbk}_{1}\leq{s}\leq{\Bbbk}_{2},
\end{equation}
\subsection{Definition }
AT of some partial derivatives is define as \cite{thermoAT,thermoAT1,thermoATF}
\begin{enumerate}
\item $\mathbf{A}\left[{{{{u}}}^{'}}({{{t}}})\right]=sU(s)-\frac{{{u}}(0)}{s}, $
\item $\mathbf{A}\left[{{{{u}}}^{''}}({{{t}}})\right]=s^{2}U(s)-\frac{{{{{u}}}^{'}}(0)}{s}-{{{u}}(0)}, $\\
\vdots
\item $\mathbf{A}\left[{{{{u}}}^{n}}({{{t}}})\right]=s^{n}U(s)-\sum_{{j}=0}^{n-1} \frac{{{{{u}}}^{k}}(0)}{s^2-n+j}, $
\end{enumerate}
where U(s) is the AT of ${{u}}{{{t}}}$.\\
 {\textbf{Remark} }\\
Let$ g$ and $j$ are two constants, ${{u}}({{{x}}},{{{t}}})$ and ${{v}}({{{x}}},{{{t}}})$ are the functions define within set A then:
\begin{equation*}\label{poo}
  \mathbf{A}\left[g{{{{u}}}}({{{x}}},{{{t}}}) \pm j{{v}}({{{x}}},{{{t}}}) \right]=g\mathbf{A}\left[{{{{u}}}}({{{x}}},{{{t}}}) \right] \pm j\mathbf{A}\left[{{v}}({{{x}}},{{{t}}})\right]
\end{equation*}
\subsection{{Definition}}\label{ap}
To control the non-linearity in a system is always challenging task for the researcher. In this connection G. Adomian and R. RACH is defined a new class of the Adomian polynomials in \cite{thermoadomian}. Which can be used to control the non-linear term in this paper, which is define as
\begin{equation}\label{s4}
N{u({{{{x}}}},{{{t}}})}=\sum_{{{j}=0}}^{\infty}A_{{j}},
\end{equation}
where
\begin{equation}
A_{{j}}=\frac{1}{{{j}}!}\left[\frac{d^{{j}}}{d\wp^{{j}}}\left\{N\sum_{{j}=0}^{\infty}(\wp^{{j}}u_{j} ({{{{x}}}},{{{t}}}))\right\}\right]_{\wp=0},
\qquad {j}=0,1,\cdots,
\end{equation}
are known as Adomian polynomials. See G. Adomian and R. RACH \cite{thermoadomian} for the the convergence of $A_{{j}}$.\\
\subsection{Definition}
The\textmd{ Mittag-Leffler} function is define as\cite{thermoML},\\
\[
E_\beta({{{{x}}}})=\sum_{{{j}=0}}^{\infty}\frac{{{{{x}}}}^{{{j}}}}{\Gamma(\beta {{j}+1})} \qquad \beta>{0}\quad {{{{x}}}}\in \mathbb{{C}},
\]
\section{ATDM Procedure }\label{s3}
In this section, we will go through the ATDM methodology for a non-linear one-dimensional thermoelasticity coupled system. The method has been recently extended to linear and non-linear fractional dynamical systems of integro-differential equations (see ref. \cite{thermothermalkhan2024application}).
\begin{equation}\label{1}
\begin{cases}
  \begin{split}
    &D_{{{{t}}}}^{{\beta+1}}{{{u}}}({{{{x}}}},{{{t}}})-a\left(\frac{\partial {{{u}}}}{\partial{{{{x}}}}},{{{v}}}\right)\frac{\partial^{2}}{\partial{{{{x}}}}^{2}}{{{u}}}({{{{x}}}},{{{t}}})+b\left(\frac{\partial {{{u}}}}{\partial{{{{x}}}}},{{{v}}}\right)\frac{\partial {{{v}}}({{{{x}}}},{{{t}}})}{\partial{{{{x}}}}}-h({{{{x}}}},{{{t}}})=0, \ \ 0<\beta\leq1,\ \ {{{{x}}}}\in \Omega,\\
    &c\left(\frac{\partial {{{u}}}}{\partial{{{{x}}}}},{{{v}}}\right)D_{{{{t}}}}^{{\beta}}{{{v}}}({{{{x}}}},{{{t}}})+\frac{\partial^{2} {{{u}}}({{{{x}}}},{{{t}}})}{\partial{{{t}}}\partial{{{{x}}}}}-d({{{v}}})\frac{\partial^{2}}{\partial{{{{x}}}}^{2}}{{{v}}}({{{{x}}}},{{{t}}})-m({{{{x}}}},{{{t}}})=0, \ \ 0<\beta\leq1,  \ \ {{{t}}}>0.
  \end{split}
  \end{cases}
\end{equation}

with the following initial conditions
\begin{equation*}\label{2}
   {{{u}}}({{{{x}}}},0)=f_{0}({{{{x}}}}),\ \ \ {{{u}}}_{{{{t}}}}({{{{x}}}},{{0}})=f_{1}({{{{x}}}}), \ \ {{{v}}}({{{{x}}}},0)=g_{0}({{{{x}}}}),
\end{equation*}
where ${{u}}({{{{x}}}},{{{t}}})$ is the displacement, ${{{v}}({{{{x}}}},{{{t}}})}$ is the temperature difference, $a\left(\frac{\partial {{{u}}}}{\partial{{{{x}}}}},{{{v}}}\right), \ c\left(\frac{\partial {{{u}}}}{\partial{{{{x}}}}},{{{v}}}\right),$ $ d({{{v}}}), \ h({{{{x}}}},{{{t}}})$ and $ \ m({{{{x}}}},{{{t}}})$ are the smooth functions.
Now let us assume the following
\begin{equation}\label{3}
\begin{split}
  &a\left(\frac{\partial {{u}}}{\partial{{{{x}}}}},{{{v}}}\right)=c\left(\frac{\partial {{{u}}}}{\partial{{{{x}}}}},{{{v}}}\right)=d({{{v}}})=1,\ \  b\left(\frac{\partial {{{u}}}}{\partial{{{{x}}}}},{{{v}}}\right)=\frac{\partial {{{u}}}}{\partial{{{{x}}}}}{{{v}}},
\end{split}
\end{equation}
putting equation (\ref{3}) into Eq. (\ref{1}), we obtain
\begin{equation}\label{4}
\begin{cases}
  \begin{split}
    &D_{{{{t}}}}^{{\beta}+1}{{{u}}}({{{{x}}}},{{{t}}})-\frac{\partial^{2}}{\partial{{{{x}}}}^{2}}{{{u}}}({{{{x}}}},{{{t}}})+{{{v}}}\frac{\partial {{{u}}}({{{{x}}}},{{{t}}})}{\partial{{{{x}}}}}\frac{\partial {{{v}}}({{{{x}}}},{{{t}}})}{\partial{{{{x}}}}}+h({{{{x}}}},{{{t}}})=0, \ \ 0<\beta\leq1, \\
    &D_{{{{t}}}}^{{\beta}}{{{v}}}({{{{x}}}},{{{t}}})+{{{v}}}({{{{x}}}},{{{t}}})\frac{\partial {{{u}}}}{\partial{{{{x}}}}}\frac{\partial^{2} {{{u}}}({{{{x}}}},{{{t}}})}{\partial{{{t}}}\partial{{{{x}}}}}-\frac{\partial^{2}}{\partial{{{{x}}}}^{2}}{{{v}}}({{{{x}}}},{{{t}}})+m({{{{x}}}},{{{t}}})=0, \ \ \ \ 0<\beta\leq1.
 \end{split}
  \end{cases}
\end{equation}
Applying AT to equation (\ref{4}), we can get
\begin{equation}\label{5}
\begin{cases}
  \begin{split}
    \mathbf{A}{{{u}}}({{{{x}}}},{{{t}}})&-\frac{f_{0}({{{{x}}}})}{s}+\frac{f_{1}({{{{x}}}})}{s^{2}}
    +\frac{1}{s^{\beta+1}}\mathbf{A}\left[-\frac{\partial^{2}}{\partial{{{{x}}}}^{2}}{{{u}}}({{{{x}}}},{{{t}}})+{{{v}}}\frac{\partial {{{u}}}({{{{x}}}},{{{t}}})}{\partial{{{{x}}}}}\frac{\partial {{{v}}}({{{{x}}}},{{{t}}})}{\partial{{{{x}}}}}\right]+\frac{1}{s^{\beta+1}}\mathbf{A}\left\{{h({{{{x}}}},{{{t}}})}\right\}=0,\\
    \mathbf{A}{{{v}}}({{{{x}}}},{{{t}}})&-\frac{g_{0}({{{{x}}}})}{s}+\frac{1}{s^{\beta}}\mathbf{A}\left[{{{v}}}({{{{x}}}},{{{t}}})\frac{\partial {{{u}}}}{\partial{{{{x}}}}}\frac{\partial^{2} {{{u}}}({{{{x}}}},{{{t}}})}{\partial{{{t}}}\partial{{{{x}}}}}-\frac{\partial^{2}}{\partial{{{{x}}}}^{2}}{{{v}}}({{{{x}}}},{{{t}}})\right]+\frac{1}{s^{\beta}}\mathbf{A}\left\{{m({{{{x}}}},{{{t}}})}\right\}=0,
  \end{split}
  \end{cases}
\end{equation}
taking inverse AT of Eq. (\ref{5}), we have
\begin{equation}\label{7}
\begin{cases}
  \begin{split}
    {{{u}}}({{{{x}}}},{{{t}}})&={f_{0}({{{{x}}}})}-{{{t}}}{f_{1}({{{{x}}}})}
    -\mathbf{A}^{-1}\left[\frac{1}{s^{\beta+1}}\mathbf{A}\left\{{h({{{{x}}}},{{{t}}})}\right\}\right]-\mathbf{A}^{-1}\left\{\frac{1}{s^{\beta+1}}\mathbf{A}\left[-\frac{\partial^{2}}{\partial{{{{x}}}}^{2}}{{{u}}}({{{{x}}}},{{{t}}})+{{{v}}}\frac{\partial {{{u}}}({{{{x}}}},{{{t}}})}{\partial{{{{x}}}}}\frac{\partial {{{v}}}({{{{x}}}},{{{t}}})}{\partial{{{{x}}}}}\right]\right\},\\
    {{{v}}}({{{{x}}}},{{{t}}})&={g_{0}({{{{x}}}})}-\mathbf{A}^{-1}\left[\frac{1}{s^{\beta}}\mathbf{A}\left\{{m({{{{x}}}},{{{t}}})}\right\}\right]-\mathbf{A}^{-1}\left\{\frac{1}{s^{\beta}}\mathbf{A}\left[{{{v}}}({{{{x}}}},{{{t}}})\frac{\partial {{{u}}}}{\partial{{{{x}}}}}\frac{\partial^{2} {{{u}}}({{{{x}}}},{{{t}}})}{\partial{{{t}}}\partial{{{{x}}}}}-\frac{\partial^{2}}{\partial{{{{x}}}}^{2}}{{{v}}}({{{{x}}}},{{{t}}})\right]\right\},
  \end{split}
  \end{cases}
\end{equation}

the  decomposition series form solution for ${{u}}({{{{x}}}},{{{t}}})$ and ${{v}}({{{{x}}}},{{{t}}})$ is given by,\\
\begin{equation}\label{eq2}
{{u}}({{{{x}}}},{{{t}}})=\sum_{{{{j}}}=0}^{\infty}{{u}}_{{{j}}}({{{{x}}}},{{{t}}}),\ \ \  \
{{v}}({{{{x}}}},{{{t}}})=\sum_{{{{j}}}=0}^{\infty}{{v}}_{{{j}}}({{{{x}}}},{{{t}}}),
\end{equation}

Non-linearities in  the given system can be represented as
\begin{equation}\label{eq3}
\begin{cases}
\begin{split}
\mathcal{A}_{{{j}}}=\frac{1}{{{{j}}}!}\left[\frac{\partial^{{{j}}}}{\partial\wp^{{{j}}}}\left\{\mathcal{N}_1\left(\sum_{k=0}^{\infty}\wp^k{{u}}_k,\sum_{k=0}^{\infty}\wp^k{{v}}_k\right)\right\}\right]_{\wp=0},\\
\mathcal{B}_{{{j}}}=\frac{1}{{{{j}}}!}\left[\frac{\partial^{{{j}}}}{\partial\wp^{{{j}}}}\left\{\mathcal{N}_2\left(\sum_{k=0}^{\infty}\wp^k{{u}}_k,\sum_{k=0}^{\infty}\wp^k{{v}}_k\right)\right\}\right]_{\wp=0}.
\end{split}
\end{cases}
\end{equation}
Using equation (\ref{eq2}) and equation (\ref{eq3}) in equation (\ref{7}), we have

\begin{equation}\label{8}
\begin{cases}
  \begin{split}
    \sum_{{{{j}}}=0}^{\infty}{{u}}_{{{j}}}({{{{x}}}},{{{t}}})&={f_{0}({{{{x}}}})}-{{{t}}}{f_{1}({{{{x}}}})}
    -\mathbf{A}^{-1}\left[\frac{1}{s^{\beta+1}}\mathbf{A}\left\{{h({{{{x}}}},{{{t}}})}\right\}\right]-\mathbf{A}^{-1}\left\{\frac{1}{s^{\beta+1}}\mathbf{A}\left[-\frac{\partial^{2}}{\partial{{{{x}}}}^{2}}\left({\sum_{{{{j}}}=0}^{\infty}{{u}}_{{{j}}}}({{{{x}}}},{{{t}}})\right)+\mathcal{A}_{{{j}}}\right]\right\},\\
    \sum_{{{{j}}}=0}^{\infty}{{v}}_{{{j}}}({{{{x}}}},{{{t}}})&={g_{0}({{{{x}}}})}-\mathbf{A}^{-1}\left[\frac{1}{s^{\beta}}\mathbf{A}\left\{{m({{{{x}}}},{{{t}}})}\right\}\right]-\mathbf{A}^{-1}\left\{\frac{1}{s^{\beta}}\mathbf{A}\left[\mathcal{B}_{{{j}}}-\frac{\partial^{2}}{\partial{{{{x}}}}^{2}}\left(\sum_{{{{j}}}=0}^{\infty}{{v}}_{{{j}}}({{{{x}}}},{{{t}}})\right)\right]\right\},
  \end{split}
  \end{cases}
\end{equation}
 the following recurrence relation can occur
 \begin{equation}\label{129}
\begin{cases}
  \begin{split}
  {{u}}_{{{0}}}({{{{x}}}},{{{t}}})&\coloneqq{f_{0}({{{{x}}}})}-{{{t}}}{f_{1}({{{{x}}}})}
    -\mathbf{A}^{-1}\left[\frac{1}{s^{\beta+1}}\mathbf{A}\left\{{h({{{{x}}}},{{{t}}})}\right\}\right],\\
    {{v}}_{{{0}}}({{{{x}}}},{{{t}}})&\coloneqq{g_{0}({{{{x}}}})}-\mathbf{A}^{-1}\left[\frac{1}{s^{\beta}}\mathbf{A}\left\{{m({{{{x}}}},{{{t}}})}\right\}\right],
  \end{split}
  \end{cases}
\end{equation}
and for ${j\geq0}$, we have
\begin{equation}\label{88}
\begin{cases}
  \begin{split}
    {{u}}_{{{j}+1}}({{{{x}}}},{{{t}}})&\coloneqq\mathbf{A}^{-1}\left\{\frac{1}{s^{\beta+1}}\mathbf{A}\left[\frac{\partial^{2}}{\partial{{{{x}}}}^{2}}\left({{{u}}_{{{j}}}}({{{{x}}}},{{{t}}})\right)-\mathcal{A}_{{{j}}}\right]\right\},\\
   {{v}}_{{{j}+1}}({{{{x}}}},{{{t}}})&\coloneqq\mathbf{A}^{-1}\left\{\frac{1}{s^{\beta}}\mathbf{A}\left[-\mathcal{B}_{{{j}}}+\frac{\partial^{2}}{\partial{{{{x}}}}^{2}}\left({{v}}_{{{j}}}({{{{x}}}},{{{t}}})\right)\right]\right\}.
  \end{split}
  \end{cases}
\end{equation}
\subsection{Example}\label{example3}
1st we consider a linear and singular thermo-elasticity coupled system of the following form
\begin{equation}\label{17}
  \begin{split}
    &\mathbf{D}_{{{t}}}^{{\beta+1}}{{{u}}}({{{x}}},{{t}})-\frac{1}{{{{x}}}^{2}}\times\frac{\partial}{\partial{{{x}}}}\left({{{x}}}^{2}\frac{\partial}{\partial{{{x}}}}{{{u}}}({{{x}}},{{t}})\right)+{{{x}}}\frac{\partial {{{v}}}({{{x}}},{{t}})}{\partial{{{x}}}}=2{{{x}}}^{2}{{t}}+{{t}}+6,\\
    &\mathbf{D}_{{{t}}}^{{\beta}}{{{v}}}({{{x}}},{{t}})-\frac{1}{{{{x}}}^{2}}\times\frac{\partial}{\partial{{{x}}}}\left({{{{{x}}}}}^2\frac{\partial}{\partial{{{x}}}}{{{v}}}({{{x}}},{{t}})\right)+\frac{\partial^{2} {{{u}}}({{{x}}},{{t}})}{\partial{{t}}\partial{{{x}}}}+6{{{x}}}^{2}-3{{{x}}}^{2}-6{{t}}=0, \ \ \ {{{t}}}>0 \ \ 0<\beta\leq1.
  \end{split}
\end{equation}
with initial conditions
\begin{equation}\label{18}
  {{{u}}}({{{x}}},0)={{{x}}}^{2},\ \ {{{u}}}_{{{t}}}({{{x}}},{{t}})={{{x}}}^{2}, \ \ {{{v}}}({{{x}}},0)=0.
\end{equation}
Recall the recurrence relation from the methodology section  \ref{s3}, particularly equation (\ref{129}) and equation (\ref{88}), we  obtained
\begin{equation*}
\begin{split}
&u_{0}({{{x}}},{{t}})\coloneqq {{{x}}}^{2}{{t}} +{{{x}}}^{2}+\frac{2 {{t}}^{\beta +1} \left(-3+\frac{{t} \left({{{x}}}^{2}-3\right)}{2+\beta}\right)}{\Gamma \! \left(2+\beta \right)},\\
&v_{0}({{{x}}},{{t}})\coloneqq 3 {{t}}^{\beta} \left(\frac{{{{x}}}^{2}}{\Gamma \! \left(\beta +1\right)}-\frac{2 t}{\Gamma \! \left(2+\beta \right)}\right),
\\
&u_{1}({{{x}}},{{t}})\coloneqq \frac{6 {{t}}^{2+\beta}}{\Gamma \! \left(3+\beta \right)}+\frac{6 {{t}}^{\beta +1}}{\Gamma \! \left(2+\beta \right)}+\frac{12 {{t}}^{3+2 \beta}}{\Gamma \! \left(4+2 \beta \right)}-\frac{6 {{{x}}}^{2} {{t}}^{1+2 \beta}}{\Gamma \! \left(2+2 \beta \right)},\\
&v_{1}({{{x}}},{{t}})\coloneqq -\frac{2 {{{x}}}^{2} {{t}}^{\beta}}{\Gamma \! \left(\beta +1\right)}-\frac{4 {{{x}}}^{2} {{t}}^{1+2 \beta}}{\Gamma \! \left(2+2 \beta \right)}+\frac{18 {{t}}^{2 \beta}}{\Gamma \! \left(1+2 \beta \right)},
\\
&u_{2}({{{x}}},{{t}})\coloneqq \frac{4 {{{x}}}^{2} {{t}}^{1+2 \beta}}{\Gamma \! \left(2+2 \beta \right)}+\frac{4 {{t}}^{2+3 \beta} \left(2 {{{x}}}^{2}-9\right)}{\Gamma \! \left(3+3 \beta \right)}\\
&v_{2}({{{x}}},{{t}})\coloneqq -\frac{12 {{t}}^{2 \beta}}{\Gamma \! \left(1+2 \beta \right)}+\frac{12 {{{x}}}^{2} {{t}}^{3 \beta}}{\Gamma \! \left(1+3 \beta \right)}-\frac{24 {{t}}^{1+3 \beta}}{\Gamma \! \left(2+3 \beta \right)},
\\
&u_{3}({{{x}}},{{t}})\coloneqq \frac{24 {{t}}^{2+3 \beta}}{\Gamma \! \left(3+3 \beta \right)}+\frac{48 {{t}}^{3+4 \beta}}{\Gamma \! \left(4+4 \beta \right)}-\frac{24 {{{x}}}^{2} {{t}}^{1+4 \beta}}{\Gamma \! \left(2+4 \beta \right)},\\
&v_{3}({{{x}}},{{t}})\coloneqq -\frac{8 {{{x}}}^{2} {{t}}^{3 \beta}}{\Gamma \! \left(1+3 \beta \right)}+\frac{72 {{t}}^{4 \beta}}{\Gamma \! \left(1+4 \beta \right)}-\frac{16 {{{x}}}^{2} {{t}}^{1+4 \beta}}{\Gamma \! \left(2+4 \beta \right)},\\
&u_{4}({{{x}}},{{t}})\coloneqq \frac{16 {{{x}}}^{2} {{t}}^{1+4 \beta}}{\Gamma \! \left(2+4 \beta \right)}+\frac{16 {{t}}^{2+5 \beta} \left(2 {{{x}}}^{2}-9\right)}{\Gamma \! \left(3+5 \beta \right)},\\
&v_{4}({{{x}}},{{t}})\coloneqq -\frac{48 {{t}}^{4 \beta}}{\Gamma \! \left(1+4 \beta \right)}+\frac{48 {{{x}}}^{2} {{t}}^{5 \beta}}{\Gamma \! \left(1+5 \beta \right)}-\frac{96 {{t}}^{1+5 \beta}}{\Gamma \! \left(2+5 \beta \right)},\\
&u_{5}({{{x}}},{{t}})\coloneqq \frac{96 {{t}}^{2+5 \beta}}{\Gamma \! \left(3+5 \beta \right)}+\frac{192 {{t}}^{3+6 \beta}}{\Gamma \! \left(4+6 \beta \right)}-\frac{96 {{{x}}}^{2} {{t}}^{1+6 \beta}}{\Gamma \! \left(2+6 \beta \right)},
\\
&v_{5}({{{x}}},{{t}})\coloneqq -\frac{32 {{{x}}}^{2} {{t}}^{5 \beta}}{\Gamma \! \left(1+5 \beta \right)}+\frac{288 {{t}}^{6 \beta}}{\Gamma \! \left(1+6 \beta \right)}-\frac{64 {{{x}}}^{2} {{t}}^{1+6 \beta}}{\Gamma \! \left(2+6 \beta \right)},\\
\end{split}
\end{equation*}
\ \ \ \ \ \ \ \ \ \ \ \ \ \ \ \ \ \ \ \ \ \ \ \ \ \ \ \ \ \ \ \ \ \ \ \ \ \ \ \ \ \ \ \ \ \ \ \ \vdots\\
\begin{equation*}
\begin{split}
&u_{100}({{{x}}},{{t}})\coloneqq \frac{1267650600228229401496703205376 {{{x}}}^{2} {{t}}^{1+100 \beta}}{\Gamma \! \left(2+100 \beta \right)}+\frac{1267650600228229401496703205376 {{t}}^{2+101 \beta} \left(2 {{{x}}}^{2}-9\right)}{\Gamma \! \left(3+101 \beta \right)},\\
&v_{100}({{{x}}},{{t}})\coloneqq -\frac{3802951800684688204490109616128 {{t}}^{100 \beta}}{\Gamma \! \left(1+100 \beta \right)}-\frac{7605903601369376408980219232256 {{t}}^{1+101 \beta}}{\Gamma \! \left(2+101 \beta \right)},\\
\end{split}
\end{equation*}
\ \ \ \ \ \ \ \ \ \ \ \ \ \ \ \ \ \ \ \ \ \ \ \ \ \ \ \ \ \ \ \ \ \ \ \ \ \ \ \ \ \ \ \ \ \ \ \ \vdots\\
It is clear that the decomposition series form solution for ${{u}}({{{{x}}}},{{{t}}})$ and ${{v}}({{{{x}}}},{{{t}}})$ is given by,\\

\begin{equation}\label{14}
\begin{split}
 {u}({{{x}}},{{t}})&=\sum _{n=0}^{\infty }{u}_{{n}} ({{{x}}},{{t}})={u}_{0}({{{x}}},{{t}})+{u}_{1}({{{x}}},{{t}})+{u}_{2}({{{x}}},{{t}})+{u}_{3}({{{x}}},{{t}})+\ldots\\
 {v}({{{x}}},{{t}})&=\sum _{n=0}^{\infty }{v}_{{n}} ({{{x}}},{{t}})={v}_{0}({{{x}}},{{t}})+{v}_{1}({{{x}}},{{t}})+{v}_{2}({{{x}}},{{t}})+{v}_{3}({{{x}}},{{t}})+\ldots\\
 \end{split}
\end{equation}
Putting $u_{0}, \  u_{1},\ \cdots,$ and $v_{0}, \ v_{1}, \ \cdots,$ in equation (\,\ref{14}\,) and setting $\beta=1$, we obtain

\begin{equation}\label{15}
\begin{split}
u({{{x}}},{t})=& {{{x}}}^{2}{{t}} +{{{x}}}^{2}+\frac{{{t}}^{2} \left({{{x}}}^{2} t -{t}3 t -{t}9\right)}{3}+\frac{{{t}}^{2} \left({{t}}^{3}-10 {{{x}}}^{2}{{t}} +10{{t}} +30\right)}{10}\\
&+\frac{{{t}}^{3} \left(2 {{t}}^{2} {{{x}}}^{2}-9 {{t}}^{2}+20 {{{x}}}^{2}\right)}{30}+\frac{{{t}}^{5} \left({{t}}^{2}-21 {{{x}}}^{2}+21\right)}{105}+\ldots~~~~\longrightarrow {{{x}}}^{2}{{t}} +{{{x}}}^{2}\\
v({{{x}}},{t})=& -3 {t} \left(-{{{x}}}^{2}+t \right)-\frac{{t} \left(2 {{t}}^{2} {{{x}}}^{2}+6 {{{x}}}^{2}-27 t \right)}{3}-{{t}}^{2} \left(-2 {{{x}}}^{2}{{t}} +{{t}}^{2}+6\right)\\
&-\frac{{{t}}^{3} \left(2 {{t}}^{2} {{{x}}}^{2}+20 {{{x}}}^{2}-45 t \right)}{15}-\frac{2 {{t}}^{4} \left(-3 {{{x}}}^{2}{{t}} +{{t}}^{2}+15\right)}{15}+\ldots~~~~\longrightarrow{{{x}}}^{2}{{t}}.
\end{split}
\end{equation}
Implies that
$u({{{x}}},{t})={{{x}}}^{2}{{t}} +{{{x}}}^{2}$ and $v({{{x}}},{t})={{{x}}}^{2}{{t}}$. Which  is an exact solution of  for problem  (\ref{example3}). In addition, we observe that the final result can be compared to the outcome achieved by HPTM, LADM, and LRPSM (see \cite{thermoNaveed,thermoShabaz,thermo32}).

\begin{figure}[H]\centering
	\includegraphics[width=8.0 cm]{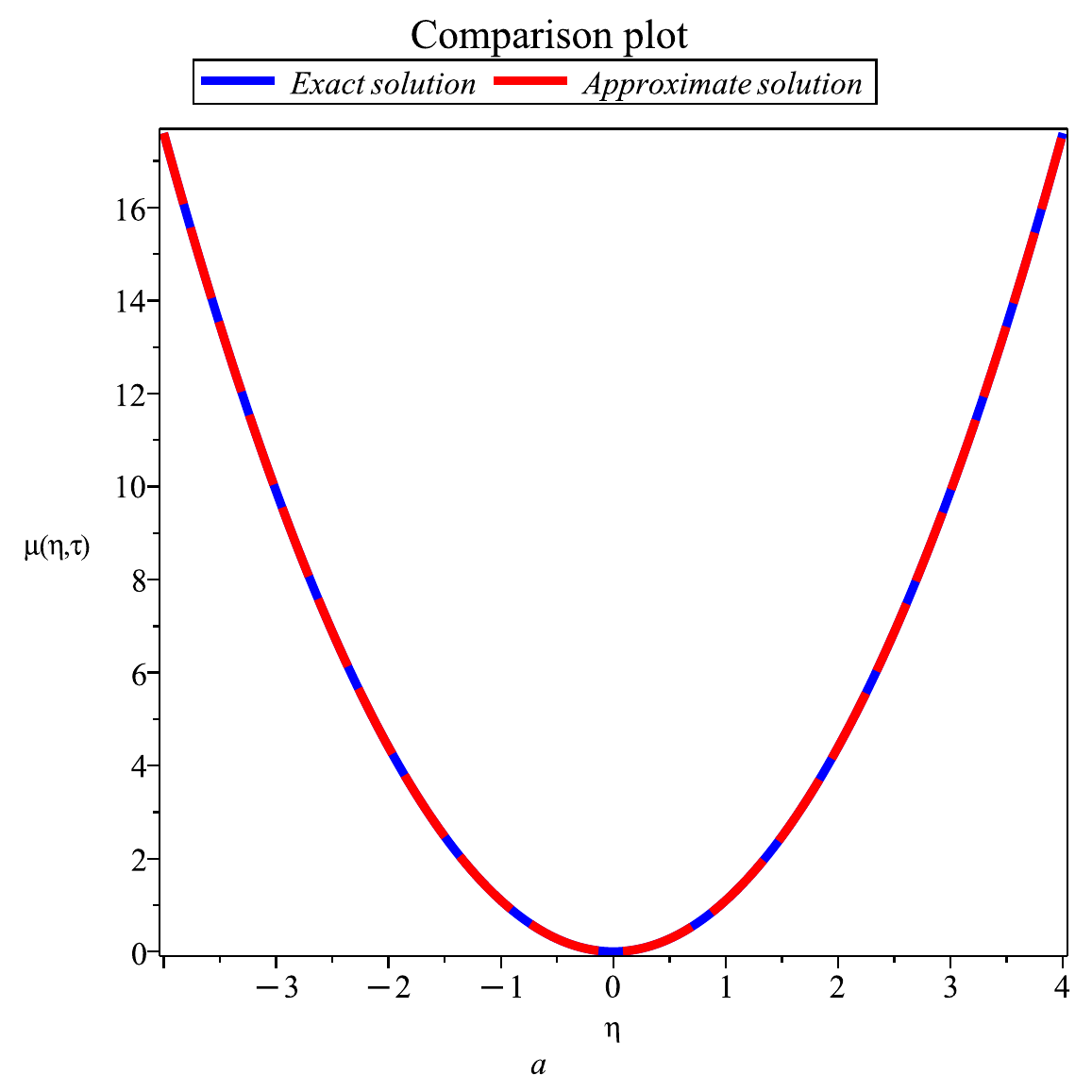}
	\includegraphics[width=8.0 cm]{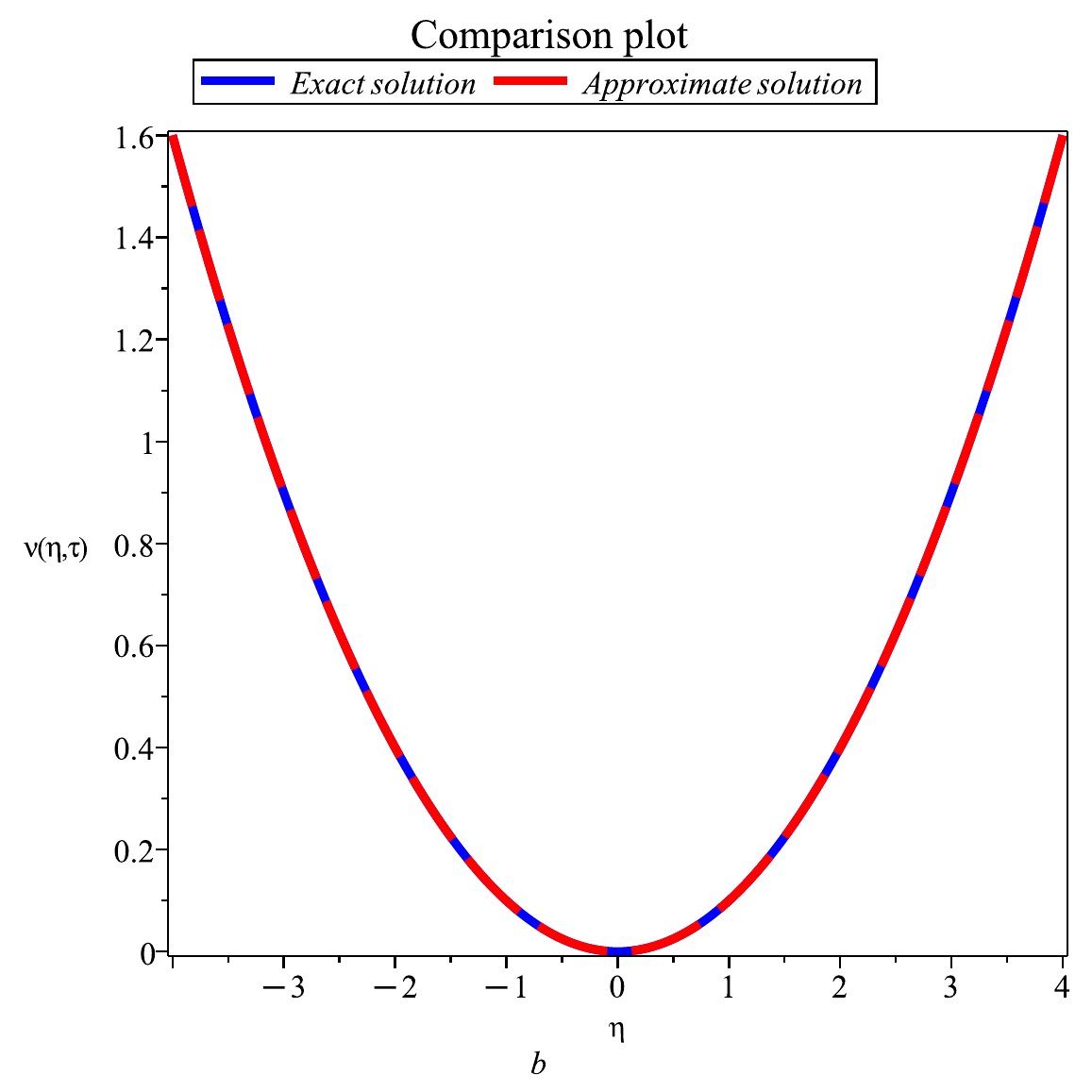}
	\caption{Comparisons  plots of approximate and exact solution for $(a)=u({{{x}}},{t})$ and $(b)=v({{{x}}},{t})$} of {example} (\ref{example3}).\label{figure3}
\end{figure}
 \begin{table}[H]
	\caption{\textbf{Comparison of {$\mathrm{L}_{2}$ Norm error}  between different iterations estimates at $\sqrt{{u}_{j+1}-{u}_{j}}$   of example (\ref{example3}).}}\label{table1}
	\centering
		\resizebox{\textwidth}{!}
	{
		\begin{tabular}{cccccc}
			\toprule
			  \textbf{Time level } 	& \textbf{$\mathrm{L}_{2}$ error}	& \textbf{$\mathrm{L}_{2}$ error}     & \textbf{$\mathrm{L}_{2}$ error}& \textbf{$\mathrm{L}_{2}$ error}&\textbf{$\mathrm{L}_{2}$ error}\\
\textbf{${{t}}$} &$\sqrt{{u}_{26}-{u}_{25}}$	& \textbf{ $\sqrt{{u}_{51}-{u}_{50}}$ }     & \textbf{  $\sqrt{{u}_{101}-{u}_{100}}$}& \textbf{  $\sqrt{{u}_{251}-{u}_{250}}$}&\textbf{ $\sqrt{{u}_{501}-{u}_{500}}$}\\
\toprule
 0.1 &$6.525\times 10^{-23}$& $3.334\times 10^{-50}$  &  $   2.019\times 10^{-114}$  &  $ 4.098\times 10^{-334}$& $3.761\times 10^{-745}$\\
 0.2 &$6.289\times 10^{-19}$  &  $1.119\times 10^{-42}$  &  $ 2.274\times 10^{-99}$  &  $  1.743\times 10^{-296}$  &  $  8.946\times 10^{-670}$\\
0.3 & $1.280\times 10^{-16}$  &  $ 2.828\times 10^{-38}$  &  $  1.450\times 10^{-90}$  &  $1.790\times 10^{-274}$  &  $  1.060\times 10^{-625}$\\
0.4 & $ 5.506\times 10^{-15}$  &  $ 2.828\times 10^{-38}$  &  $  2.561\times 10^{-84}$  &  $ 7.416\times 10^{-259}$  &  $  1.895\times 10^{-594}$\\
 0.5 &$1.015\times 10^{-13}$  &  $ 2.828\times 10^{-38}$  &  $  1.795\times 10^{-79}$  &  $ 9.638\times 10^{-247}$  &  $ 3.141\times 10^{-570}$\\
0.6 & $1.095\times 10^{-12}$  &  $9.956\times 10^{-33}$  &  $1.634\times 10^{-75}$  &  $ 7.615\times 10^{-237}$  &  $ 1.800\times 10^{-550}$\\
0.7 & $8.176\times 10^{-12}$  &  $ 9.503\times 10^{-31}$  &  $ 3.637\times 10^{-72}$  &  $ 1.779\times 10^{-228}$  &  $  8.050\times 10^{-534}$\\
0.8 &$ 4.662\times 10^{-11}$  &  $4.485\times 10^{-29}$  &  $2.887\times 10^{-69}$  &  $ 3.156\times 10^{-221}$  &  $ 1.393\times 10^{-519}$\\
0.9 &$ 2.165\times 10^{-10}$&$1.264\times 10^{-27}$  &  $ 1.043\times 10^{-66}$  &  $  7.820\times 10^{-215}$  &  $  1.210\times 10^{-506}$\\
1.0 & $ 8.545\times 10^{-10}$&$3.403\times 10^{-26}$  &  $ 2.024\times 10^{-64}$&$ 4.102\times 10^{-209}$&$5.617\times 10^{-495}$\\
\bottomrule
\end{tabular}}
\end{table}
\begin{figure}[H]\centering
	\includegraphics[width=8.0 cm]{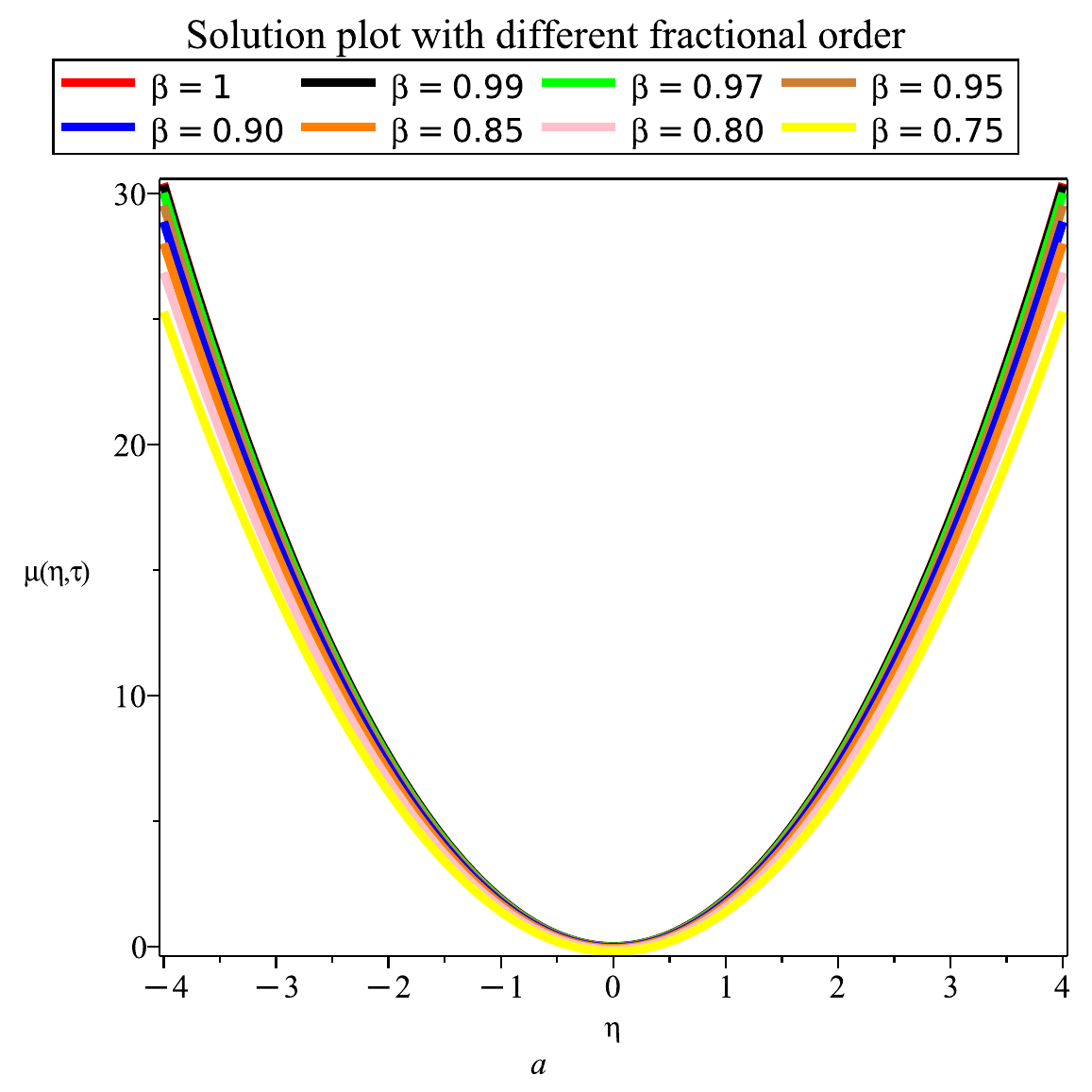}
	\includegraphics[width=8.0 cm]{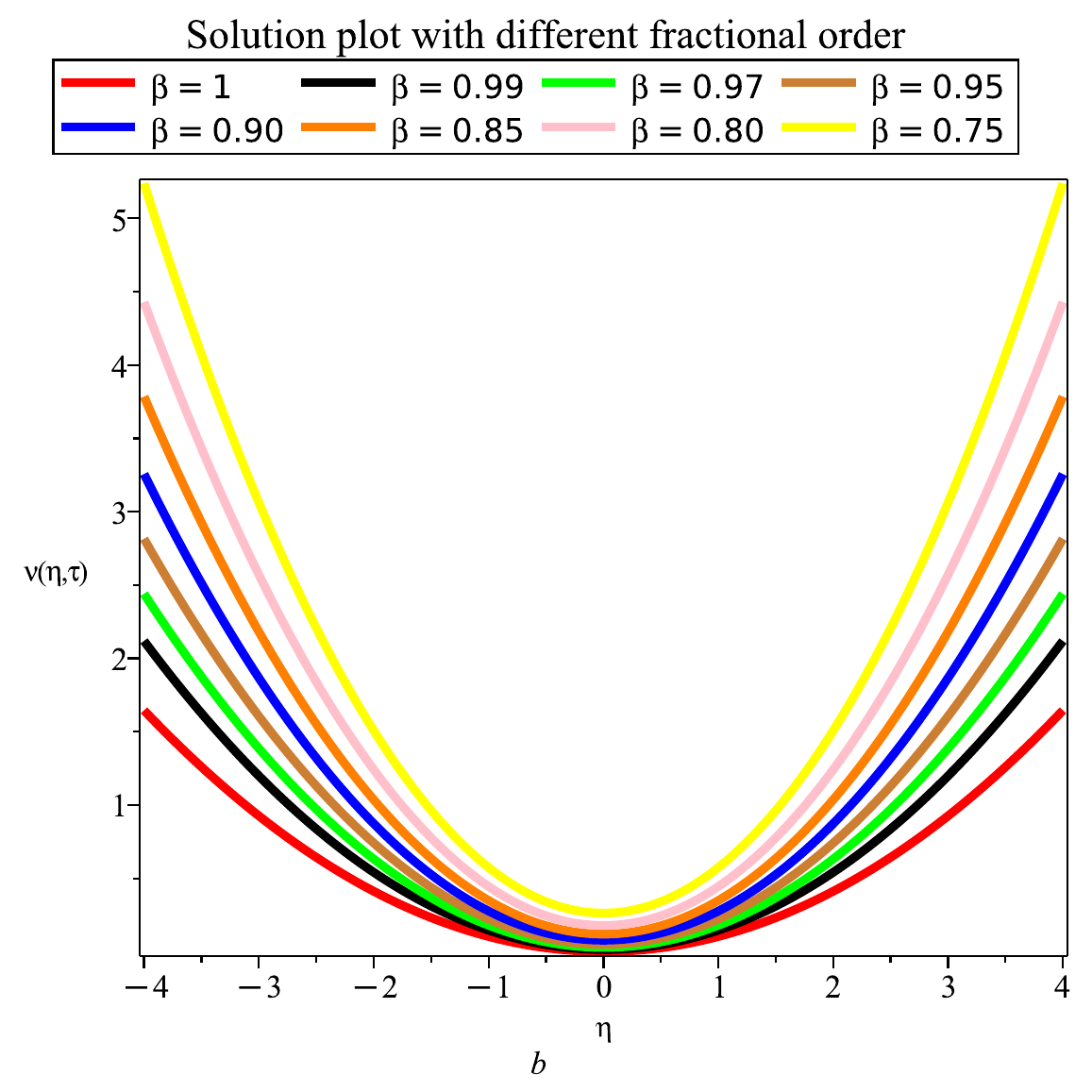}\\
\includegraphics[width=8.0 cm]{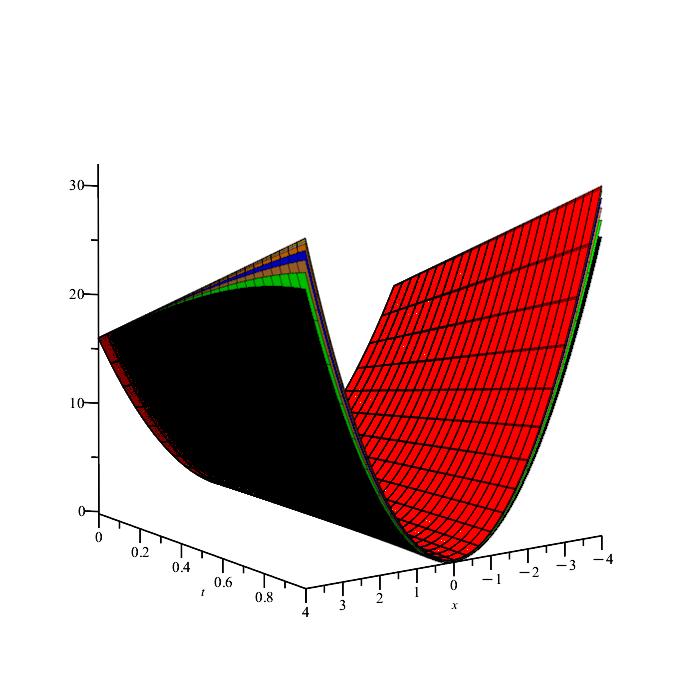}
\includegraphics[width=8.0 cm]{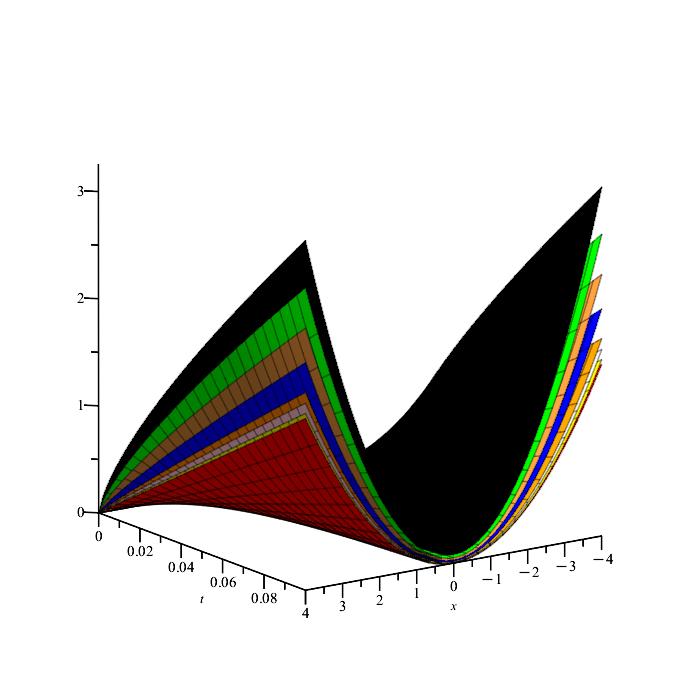}\\
\includegraphics[width=8.0 cm]{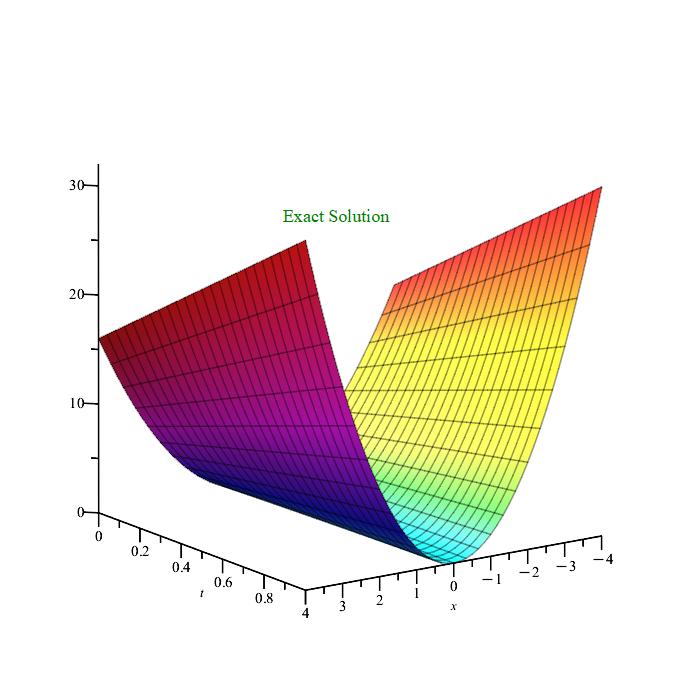}
\includegraphics[width=8.0 cm]{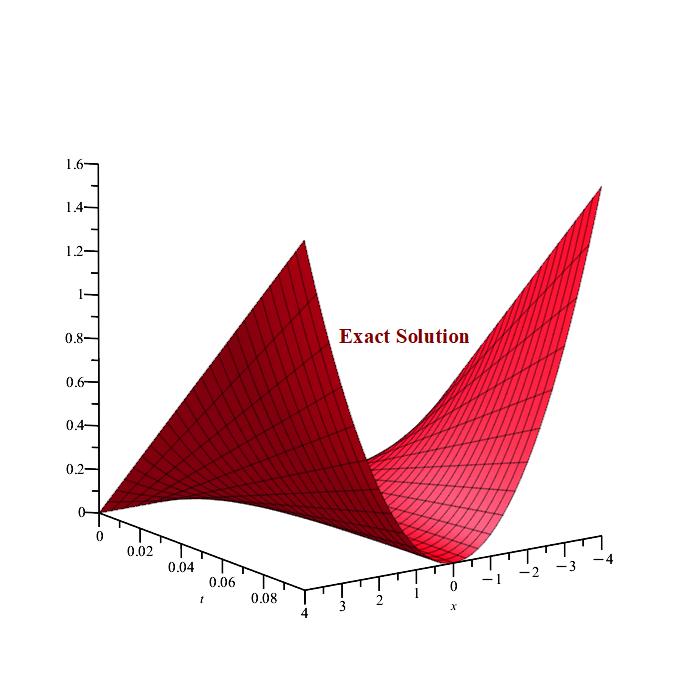}\\
	\caption{Comparisons 2D and 3D plotting for $u({{{x}}},{t})$ and $v({{{x}}},{t})$ with fractional order $\beta$ of  {example} (\ref{example3}).}\label{figure3}
\end{figure}

\begin{table}[H]
	\caption{\textbf{Numerical comparison of exact and approximate solutions with absolute error (AE) for various   time levels ${{t}}$  and different fractional  orders $\beta$ of example (\ref{example3}).} }\label{table2}
	\centering
	\resizebox{\textwidth}{!}
	{
		\begin{tabular}{cccccccc}
			\toprule
			\textbf{{Variable}} & \textbf{time level} & \textbf{Exact } 	& \textbf{Approximate }	& \textbf{Absolute error}     & \textbf{Absolute error}& \textbf{Absolute error}&\textbf{Absolute error}\\
			\textbf{${u}$, ${v}$} & \textbf{${{t}}$} & \textbf{ solutions} 	& \textbf{solutions }	& \textbf{ $\beta=0.97$ }     & \textbf{ at $\beta=0.98$}& \textbf{ at $ \beta=0.99$}&\textbf{at $\beta=1$}\\
			\midrule  \midrule
 &0.01 & 16.1600000000 & 16.1600000000   &0.0000012200   &0.0000007400   &0.0000075500   &0.0000000000\\
 &0.02 & 16.3200000000 & 16.3200000000   &0.0000083300   &0.0000051500   &0.0000478200   &0.0000000000\\
 &0.03 & 16.4799999900 & 16.4800000000   &0.0000254400   &0.0000158100   &0.0001397100   &0.0000000100\\
 &0.04 & 16.6400000000 & 16.6400000000   &0.0000559400   &0.0000349000   &0.0002976300   &0.0000000000\\
 &0.05 & 16.8000000000 & 16.8000000000   &0.0001027800   &0.0000643300   &0.0005337200   &0.0000000000\\
 &0.06 & 16.9600000000 & 16.9600000000   &0.0001686100   &0.0001058400   &0.0008585300   &0.0000000000\\
 &0.07  & 17.1199999900  & 17.1200000000   &0.0002559000   &0.0001610000   &0.0012814200   &0.0000000100\\
 &0.08  & 17.2800000000  & 17.2800000000   &0.0003668600   &0.0002312800   &0.0018109100   &0.0000000000\\
 &0.09  & 17.4400000000  & 17.4400000000   &0.0005036000   &0.0003180600   &0.0024548000   &0.0000000000\\
 &0.10  & 17.6000000000  & 17.6000000000   &0.0006680700   &0.0004226100   &0.0032202800   &0.0000000000\\
 &0.11  & 17.7600000000  & 17.7600000000   &0.0008621600   &0.0005461500   &0.0041141000   &0.0000000000\\
 &0.12  & 17.9199999900  & 17.9200000000   &0.0010875700   &0.0006898400   &0.0051425800   &0.0000000100\\
 &0.13  & 18.0800000000  & 18.0800000000   &0.0013459900   &0.0008547900   &0.0063116800   &0.0000000000\\
 &0.14  & 18.2400000000  & 18.2400000000   &0.0016390600   &0.0010420400   &0.0076270800   &0.0000000000\\
 &0.15  & 18.4000000100  & 18.4000000000   &0.0019682700   &0.0012526100   &0.0090941700   &0.0000000100\\
 &0.16  & 18.5600000000  & 18.5600000000   &0.0023350700   &0.0014874400   &0.0107181400   &0.0000000000\\
 &0.17  & 18.7200000000  & 18.7200000000   &0.0027409000   &0.0017475000   &0.0125039600   &0.0000000000\\
$\textbf{u}({{{x}}},{t})$ &0.18  & 18.8800000000  & 18.8800000000   &0.0031871500   &0.0020336800   &0.0144564500   &0.0000000000\\
 &0.19  &19.0400000000  &19.0400000000   &0.0036751000   &0.0023468600   &0.0165802800   &0.0000000000\\
 &0.20  &19.2000000000  &19.2000000000   &0.0042060700   &0.0026879000   &0.0188799500   &0.0000000000\\
 &0.21  &19.3600000000  &19.3600000000   &0.0047812800   &0.0030575900   &0.0213599000   &0.0000000000\\
 &0.22  &19.5200000000  &19.5200000000   &0.0054020000   &0.0034567900   &0.0240244200   &0.0000000000\\
 &0.23  &19.6800000100  &19.6800000000   &0.0060693800   &0.0038862600   &0.0268777800   &0.0000000100\\
 &0.24  &19.8400000100  &19.8400000000   &0.0067846100   &0.0043467900   &0.0299241200   &0.0000000100\\
 &0.25  &20.0000000100  &20.0000000000   &0.0075488000   &0.0048390800   &0.0331675000   &0.0000000100\\
 &0.26  &20.1600000200  &20.1600000000   &0.0083631300   &0.0053639600   &0.0366119900   &0.0000000200\\
 &0.27  &20.3200000400  &20.3200000000   &0.0092286600   &0.0059221000   &0.0402616200   &0.0000000400\\
 &0.28  &20.4800000400  &20.4800000000   &0.0101465200   &0.0065142300   &0.0441203300   &0.0000000400\\
 &0.29  &20.6400000600  &20.6400000000   &0.0111177500   &0.0071410800   &0.0481920400   &0.0000000600\\
 &0.30  &20.8000000800  &20.8000000000   &0.0121434400   &0.0078033400   &0.0524807300   &0.0000000800\\
 &0.31  &20.9600001100  &20.9600000000   &0.0132246400   &0.0085017100   &0.0569902700   &0.0000001100\\
 &0.32 & 21.1200001600 & 21.1200000000   &0.0143623900   &0.0092368600   &0.0617246000   &0.0000001600\\
 &0.33 & 21.2800002200 & 21.2800000000   &0.0155577500   &0.0100095000   &0.0666875900   &0.0000002200\\
 &0.34 & 21.4400002800 & 21.4400000000   &0.0168117400   &0.0108203000   &0.0718832000   &0.0000002800\\
 &0.35 & 21.6000003500 & 21.6000000000   &0.0181254000   &0.0116699300   &0.0773153400   &0.0000003500\\
 \bottomrule
 &0.01  &0.1600000000   &0.1600000000   &0.0260616574   &0.0169440742   &0.0082636890   &0.0000000000\\
 &0.02  &0.3200000000   &0.3200000000   &0.0445698858   &0.0290809413   &0.0142329553   &0.0000000000\\
 &0.03  &0.4800000000   &0.4800000000   &0.0603898407   &0.0394835526   &0.0193632857   &0.0000000000\\
 &0.04  &0.6400000001   &0.6400000000   &0.0745175863   &0.0487890496   &0.0239602596   &0.0000000001\\
 &0.05  &0.8000000003   &0.8000000000   &0.0874214477   &0.0572982191   &0.0281687288   &0.0000000003\\
 &0.06  &0.9600000011   &0.9600000000   &0.0993758710   &0.0651878280   &0.0320740543   &0.0000000011\\
 &0.07  &1.1200000030   &1.1200000000   &0.1105618040   &0.0725746730   &0.0357327730   &0.0000000030\\
 &0.08  &1.2800000090   &1.2800000000   &0.1211078300   &0.0795419120   &0.0391851920   &0.0000000090\\
 &0.09  &1.4400000190   &1.4400000000   &0.1311102080   &0.0861518890   &0.0424616010   &0.0000000190\\
 &0.10  &1.6000000410   &1.6000000000   &0.1406438410   &0.0924531930   &0.0455856190   &0.0000000410\\
 &0.11  &1.7600000790   &1.7600000000   &0.1497687930   &0.0984848330   &0.0485762530   &0.0000000790\\
 &0.12  &1.9200001450   &1.9200000000   &0.1585343850   &0.1042788740   &0.0514491360   &0.0000001450\\
 &0.13  &2.0800002540   &2.0800000000   &0.1669818820   &0.1098621900   &0.0542173900   &0.0000002540\\
 &0.14  &2.2400004280   &2.2400000000   &0.1751464390   &0.1152576870   &0.0568922150   &0.0000004280\\
&0.15  &2.4000006920   &2.4000000000   &0.1830583560   &0.1204851300   &0.0594832920   &0.0000006920\\
 &0.16  &2.5600010870   &2.5600000000   &0.1907441290   &0.1255618260   &0.0619990950   &0.0000010870\\
 &0.17  &2.7200016600   &2.7200000000   &0.1982271110   &0.1305030580   &0.0644471480   &0.0000016600\\
 \textbf{$v({{{x}}},{t})$} &0.18  &2.8800024780   &2.8800000000   &0.2055281470   &0.1353224410   &0.0668341490   &0.0000024780\\
 &0.19  &3.0400036160   &3.0400000000   &0.2126659350   &0.1400322630   &0.0691661660   &0.0000036160\\
 &0.20  &3.2000051770   &3.2000000000   &0.2196574500   &0.1446436350   &0.0714486980   &0.0000051770\\
 &0.21  &3.3600072830   &3.3600000000   &0.2265181360   &0.1491667380   &0.0736868110   &0.0000072830\\
 &0.22  &3.5200100840   &3.5200000000   &0.2332622230   &0.1536109250   &0.0758851570   &0.0000100840\\
 &0.23  &3.6800137600   &3.6800000000   &0.2399028480   &0.1579848660   &0.0780480940   &0.0000137600\\
 &0.24  &3.8400185320   &3.8400000000   &0.2464522320   &0.1622966440   &0.0801797190   &0.0000185320\\
 &0.25  &4.0000246560   &4.0000000000   &0.2529218420   &0.1665538620   &0.0822838750   &0.0000246560\\
 &0.26  &4.1600324380   &4.1600000000   &0.2593224310   &0.1707636690   &0.0843642530   &0.0000324380\\
 &0.27  &4.3200422360   &4.3200000000   &0.2656642520   &0.1749329080   &0.0864243760   &0.0000422360\\
 &0.28  &4.4800544690   &4.4800000000   &0.2719570310   &0.1790680620   &0.0884676490   &0.0000544690\\
 &0.29  &4.6400696190   &4.6400000000   &0.2782100690   &0.1831754460   &0.0904974460   &0.0000696190\\
 &0.30  &4.8000882440   &4.8000000000   &0.2844323920   &0.1872611190   &0.0925169760   &0.0000882440\\
 &0.31 & 4.9601109850 &  4.9600000000   &0.2906327300   &0.1913309920   &0.0945295080   &0.0001109850\\
 &0.32  &5.1201385740   &5.1200000000   &0.2968195390   &0.1953909140   &0.0965382390   &0.0001385740\\
 &0.33  &5.2801718410   &5.2800000000   &0.3030012040   &0.1994465590   &0.0985464220   &0.0001718410\\
 &0.34  &5.4402117270   &5.4400000000   &0.3091859080   &0.2035036780   &0.1005573230   &0.0002117270\\
 &0.35  &5.6002592970   &5.6000000000   &0.3153818360   &0.2075679110   &0.1025742510   &0.0002592970\\
 \bottomrule
	\end{tabular}}
\end{table}

\subsection{Example}\label{example32}
Consider the non-linear coupled system of the thermo-elasticity
\begin{equation}\label{14}
\begin{cases}
  \begin{split}
    &D_{{{{t}}}}^{{\beta+1}}{{{u}}}({{{{x}}}},{{{t}}})-\frac{\partial}{\partial{{{{x}}}}}\left({{{v}}}({{{{x}}}},{{{t}}})\frac{\partial}{\partial{{{{x}}}}}{{{u}}}({{{{x}}}},{{{t}}})\right)+\frac{\partial {{{v}}}({{{{x}}}},{{{t}}})}{\partial{{{{x}}}}}-2{{{{x}}}}-6{{{{x}}}}^{2}-2{{{t}}}^{2}-2=0,\ \  \ \ 0<\beta\leq1,\\
    &D_{{{{t}}}}^{{\beta}}{{{v}}}({{{{x}}}},{{{t}}})+\frac{\partial}{\partial{{{{x}}}}}\left({u}
    ({{{{x}}}},{{{t}}})\frac{\partial}{\partial{{{{x}}}}}{{{v}}}({{{{x}}}},{{{t}}})\right)+\frac{\partial^{2} {{{u}}}({{{{x}}}},{{{t}}})}{\partial{{{t}}}\partial{{{{x}}}}}+6{{{{x}}}}^{2}-2{{{t}}}^{2}-2{{{t}}}=0, \ \ 0<\beta\leq1,  \ \ {{{t}}}>0.
  \end{split}
  \end{cases}
\end{equation}

With initial conditions
\begin{equation}\label{15}
  {{{u}}}({{{{x}}}},0)={{{{x}}}}^{2},\ \ {{{u}}}_{{{{t}}}}({{{{x}}}},{{{t}}})=0, \ \ {{{v}}}({{{{x}}}},0)={{{{x}}}}^{2}.
\end{equation}
In system (\ref{14}), we have two non-linear terms ${{{v}}}({{{{x}}}},{{{t}}})\frac{\partial}{\partial{{{{x}}}}}{{{u}}}({{{{x}}}},{{{t}}})$ and ${{{u}}}({{{{x}}}},{{{t}}})\frac{\partial}{\partial{{{{x}}}}}{{{v}}}({{{{x}}}},{{{t}}})$. Let ${\mathbf{N}_1({u},{v})}$ and $\mathbf{N}_2({u},{v})$ be the two non-linear terms. Now recall the definition \ref{ap} of the Adomian polynomial to control the non linearity i.e  ${\mathbf{N}_1({u},{v})}$ and ${\mathbf{N}_2({u},{v})}$.\\
Where
\begin{equation}
{\mathbf{N}_1({u},{v})}=\sum_{{{j}=0}}^{\infty}A_{{j}},\ \ \
{\mathbf{N}_2({u},{v})}=\sum_{{{j}=0}}^{\infty}B_{{j}}.
\end{equation}
\begin{equation*}
  \begin{split}
&A_{0}\coloneqq v_{0}\! \left({{{x}}} ,{t} \right) \left(\frac{\partial}{\partial x}u_{0}\! \left({{{x}}} ,{t} \right)\right),
\\& B_{0}\coloneqq u_{0}\! \left({{{x}}} ,{t} \right) \left(\frac{\partial}{\partial x}v_{0}\! \left({{{x}}} ,{t} \right)\right),
\\&A_{1}\coloneqq v_{1}\! \left({{{x}}} ,{t} \right) \left(\frac{\partial}{\partial x}u_{0}\! \left({{{x}}} ,{t} \right)\right)+v_{0}\! \left({{{x}}} ,{t} \right) \left(\frac{\partial}{\partial x}u_{1}\! \left({{{x}}} ,{t} \right)\right),
\\&B_{1}\coloneqq u_{1}\! \left({{{x}}} ,{t} \right) \left(\frac{\partial}{\partial x}v_{0}\! \left({{{x}}} ,{t} \right)\right)+u_{0}\! \left({{{x}}} ,{t} \right) \left(\frac{\partial}{\partial x}v_{1}\! \left({{{x}}} ,{t} \right)\right),
\\&A_{2}\coloneqq v_{2}\! \left({{{x}}} ,{t} \right) \left(\frac{\partial}{\partial x}u_{0}\! \left({{{x}}} ,{t} \right)\right)+v_{1}\! \left({{{x}}} ,{t} \right) \left(\frac{\partial}{\partial x}u_{1}\! \left({{{x}}} ,{t} \right)\right)+v_{0}\! \left({{{x}}} ,{t} \right) \left(\frac{\partial}{\partial x}u_{2}\! \left({{{x}}} ,{t} \right)\right)
,\\&B_{2}\coloneqq u_{2}\! \left({{{x}}} ,{t} \right) \left(\frac{\partial}{\partial x}v_{0}\! \left({{{x}}} ,{t} \right)\right)+u_{1}\! \left({{{x}}} ,{t} \right) \left(\frac{\partial}{\partial x}v_{1}\! \left({{{x}}} ,{t} \right)\right)+u_{0}\! \left({{{x}}} ,{t} \right) \left(\frac{\partial}{\partial x}v_{2}\! \left({{{x}}} ,{t} \right)\right)
,\\
\vdots
\\&A_{n}\coloneqq v_{n}\! \left({{{x}}} ,{t} \right) \left(\frac{\partial}{\partial x}u_{0}\! \left({{{x}}} ,{t} \right)\right)+v_{n-1}\! \left({{{x}}} ,{t} \right) \left(\frac{\partial}{\partial x}u_{1}\! \left({{{x}}} ,{t} \right)\right)+\cdots+v_{0}\! \left({{{x}}} ,{t} \right) \left(\frac{\partial}{\partial x}u_{n}\! \left({{{x}}} ,{t} \right)\right)
,\\&B_{n}\coloneqq u_{n}\! \left({{{x}}} ,{t} \right) \left(\frac{\partial}{\partial x}v_{0}\! \left({{{x}}} ,{t} \right)\right)+u_{n-1}\! \left({{{x}}} ,{t} \right) \left(\frac{\partial}{\partial x}v_{1}\! \left({{{x}}} ,{t} \right)\right)+\cdots+u_{0}\! \left({{{x}}} ,{t} \right) \left(\frac{\partial}{\partial x}v_{n}\! \left({{{x}}} ,{t} \right)\right), \forall \ \ {n} \ \ \in \ \ \mathbb{N}^{+}
  \end{split}
\end{equation*}
Now, using the recurrence relation from the methodology section  \ref{s3}, particularly equation (\ref{129}) and equation (\ref{88}), we  obtained the following

\begin{equation*}
 \begin{split}
&u_{0}({{{x}}},{{t}})\coloneqq {{{x}}}^{2}
,\\&v_{0}({{{x}}},{{t}})\coloneqq {{{x}}}^{2}
,\\&u_{1}({{{x}}},{{t}})\coloneqq \frac{\left(-2 \beta^{2}-4 {{t}}^{2}-10 \beta -12\right) {{t}}^{\beta +1}}{\Gamma \! \left(4+\beta \right)}
,\\&v_{1}({{{x}}},{{t}})\coloneqq \frac{\left(2 \beta +4{{t}} +4\right) {{t}}^{\beta +1}}{\Gamma \! \left(3+\beta \right)}
,\\&v_{2}({{{x}}},{{t}})\coloneqq -\frac{8 {{t}}^{1+2 \beta} \left(2 \beta^{2}+{{t}}^{2}+5 \beta +3\right)}{\Gamma \! \left(4+2 \beta \right)}
,\\&u_{3}({{{x}}},{{t}})\coloneqq \frac{\left(-72 \beta^{2}-16 {{t}}^{2}-168 \beta -96\right) {{t}}^{2+3 \beta}}{\Gamma \! \left(5+3 \beta \right)},\\
&v_{3}({{{x}}},{{t}})\coloneqq \frac{\left(24+24 \beta +16 t \right) {{t}}^{2+3 \beta}}{\Gamma \! \left(4+3 \beta \right)}
,
\end{split}
\end{equation*}
similarly we can get
\begin{equation*}
 \begin{split}
&u_{4}({{{x}}},{{t}})\coloneqq \frac{32 {{t}}^{3+4 \beta} \left(2+2 \beta +t \right)}{\Gamma \! \left(4 \beta +5\right)}
,\\&v_{4}({{{x}}},{{t}})\coloneqq -\frac{32 {{t}}^{2+4 \beta} \left(8 \beta^{2}+{{t}}^{2}+14 \beta +6\right)}{\Gamma \! \left(4 \beta +5\right)}\\
\vdots
\end{split}
\end{equation*}
It is clear that the decomposition series form solution for ${{u}}({{{{x}}}},{{{t}}})$ and ${{v}}({{{{x}}}},{{{t}}})$ is given by,
\begin{equation}\label{1q4}
\begin{split}
 & {u}({{{x}}},{{t}})=\sum _{n=0}^{\infty }{u}_{{n}} ({{{x}}},{{t}})={u}_{0}({{{x}}},{{t}})+{u}_{1}({{{x}}},{{t}})+{u}_{2}({{{x}}},{{t}})+{u}_{3}({{{x}}},{{t}})+\ldots\\
 &{v}({{{x}}},{{t}})=\sum _{n=0}^{\infty }{v}_{{n}} ({{{x}}},{{t}})={v}_{0}({{{x}}},{{t}})+{v}_{1}({{{x}}},{{t}})+{v}_{2}({{{x}}},{{t}})+{v}_{3}({{{x}}},{{t}})+\ldots\\
 \end{split}
\end{equation}
Putting $u_{0}, \  u_{1},\ \cdots,$ and $v_{0}, \ v_{1}, \ \cdots,$ in equation (\,\ref{1q4}\,) and setting $\beta=1$, we obtain

\begin{equation}\label{1q5}
\begin{split}
u({{{x}}},{t})=&{{{x}}}^{2}-{{t}}^{2}+\frac{{{t}}^{5}}{15}-\frac{{{t}}^{5} \left({{t}}^{2}+21\right)}{315}+\frac{{{t}}^{7} \left(4+{t} \right)}{1260}+\ldots~~~~~~~~~~~~~~~~~~ ~~\longrightarrow {{{{x}}}}^{2}-{{{t}}}^{2},\\
v({{{x}}},{t})=&  {{{x}}}^{2}+\frac{2 {{t}}^{3}}{3}+{{t}}^{2}-\frac{{{t}}^{3} \left({{t}}^{2}+10\right)}{15}+\frac{{{t}}^{5} \left(t +3\right)}{45}-\frac{{{t}}^{6} \left({{t}}^{2}+28\right)}{1260}
+\ldots~ \longrightarrow {{{{x}}}}^{2}+{{{t}}}^{2}.
\end{split}
\end{equation}
Implies that
${{{u}}}({{{{x}}}},{{{t}}})={{{{x}}}}^{2}-{{{t}}}^{2}$ and ${{{v}}}({{{{x}}}},{{{t}}})={{{{x}}}}^{2}+{{{t}}}^{2}$. Which  is an exact solution of  for problem (\ref{14}). In addition, we observe that the final result can be compared to the outcome achieved by HPTM, LADM, and LRPSM (see \cite{thermoNaveed,thermoShabaz,thermo32}).
\begin{table}[H]
	\caption{\textbf{Numerical comparison of absolute error with LRPSM \cite{thermoShabaz} at ${{t}=0.1}$  different for fractional  orders $\beta$ of example (\ref{example32}).}}\label{table3}
	\centering
		\resizebox{\textwidth}{!}
	{
		\begin{tabular}{cccccccc}
			\toprule
			\textbf{{Variable}} & \textbf{Space } &\textbf{Proposed method }&\textbf{ LRPSM \cite{thermoShabaz}}& \textbf{Proposed method }&\textbf{ LRPSM \cite{thermoShabaz}}& \textbf{Proposed method }&\textbf{ LRPSM \cite{thermoShabaz}}\\
			\textbf{${u}$, ${v}$} & \textbf{${x}$} & \textbf{AE at $ \beta=0.5$}&\textbf{AE at $\beta=0.5$}&\textbf{AE at $ \beta=0.7$}&\textbf{AE at $\beta=0.7$}& \textbf{AE at $ \beta=0.9$}&\textbf{AE at $\beta=0.9$}\\
			\midrule  \midrule
 &0.0  &0.03720168090  &0.03757664310  &0.01574407768  &0.01583389134  &0.00376820694  &0.00377866212\\
&0.1  &0.03720168090  &0.03757664310  &0.01574407768  &0.01583389134  &0.00376820694  &0.00377866212\\
&0.2  &0.03720168090  &0.03757664310  &0.01574407768  &0.01583389134  &0.00376820694  &0.00377866212\\
&0.3  &0.03720168090  &0.03757664310  &0.01574407768  &0.01583389134  &0.00376820694  &0.00377866212\\
&0.4  &0.03720168090  &0.03757664310  &0.01574407760  &0.01583389130  &0.00376820700  &0.00377866210\\
\textbf{$u({{{x}}},0.1)$}&0.5  &0.03720168090  &0.03757664310  &0.01574407760  &0.01583389130  &0.00376820700  &0.00377866210\\
&0.6  &0.03720168090  &0.03757664310  &0.01574407760  &0.01583389130  &0.00376820700  &0.00377866210\\
&0.7  &0.03720168090  &0.03757664310  &0.01574407760  &0.01583389130  &0.00376820700  &0.00377866210\\
&0.8  &0.03720168090  &0.03757664310  &0.01574407760  &0.01583389130  &0.00376820700  &0.00377866210\\
&0.9  &0.03720168090  &0.03757664310  &0.01574407760  &0.01583389130  &0.00376820700  &0.00377866210\\
&1.0  &0.03720168090  &0.03757664310  &0.01574407760  &0.01583389130  &0.00376820700  &0.00377866210\\
\midrule
&0.0  &0.02150988968  &0.03757664310  &0.01241949343  &0.01583389134  &0.00337907599  &0.00377866212\\
&0.1  &0.02150988968  &0.03757664310  &0.01241949343  &0.01583389134  &0.00337907599  &0.00377866212\\
&0.2  &0.02150988968  &0.03757664310  &0.01241949343  &0.01583389134  &0.00337907599  &0.00377866212\\
&0.3  &0.02150988970  &0.03757664310  &0.01241949340  &0.01583389134  &0.00337907610  &0.00377866212\\
&0.4  &0.02150988970  &0.03757664310  &0.01241949340  &0.01583389130  &0.00337907610  &0.00377866210\\
\textbf{$v({{{x}}},0.1)$}&0.5  &0.02150988970  &0.03757664310  &0.01241949340  &0.01583389130  &0.00337907610  &0.00377866210\\
&0.6  &0.02150988970  &0.03757664310  &0.01241949340  &0.01583389130  &0.00337907610  &0.00377866210\\
&0.7  &0.02150988970  &0.03757664310  &0.01241949340  &0.01583389130  &0.00337907610  &0.00377866210\\
&0.8  &0.02150988970  &0.03757664310  &0.01241949340  &0.01583389130  &0.00337907610  &0.00377866210\\
&0.9  &0.02150988970  &0.03757664310  &0.01241949340  &0.01583389130  &0.00337907610  &0.00377866210\\
&1.0  &0.02150989100  &0.03757664310  &0.01241949400  &0.01583389130  &0.00337907600  &0.00377866210\\
 \bottomrule
	\end{tabular}}
\end{table}
\begin{figure}[H]\centering
	\includegraphics[width=6.0 cm]{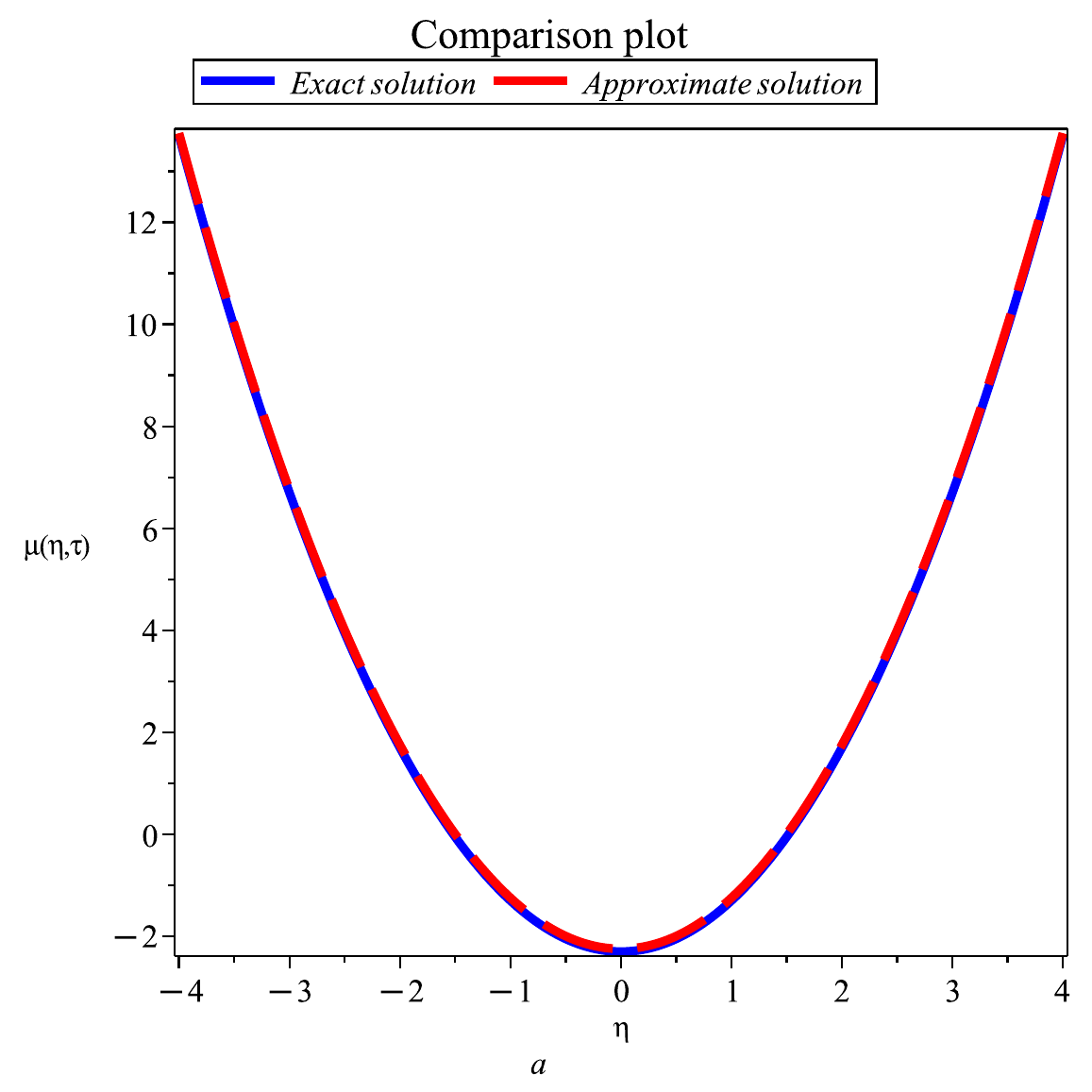}
	\includegraphics[width=6.0 cm]{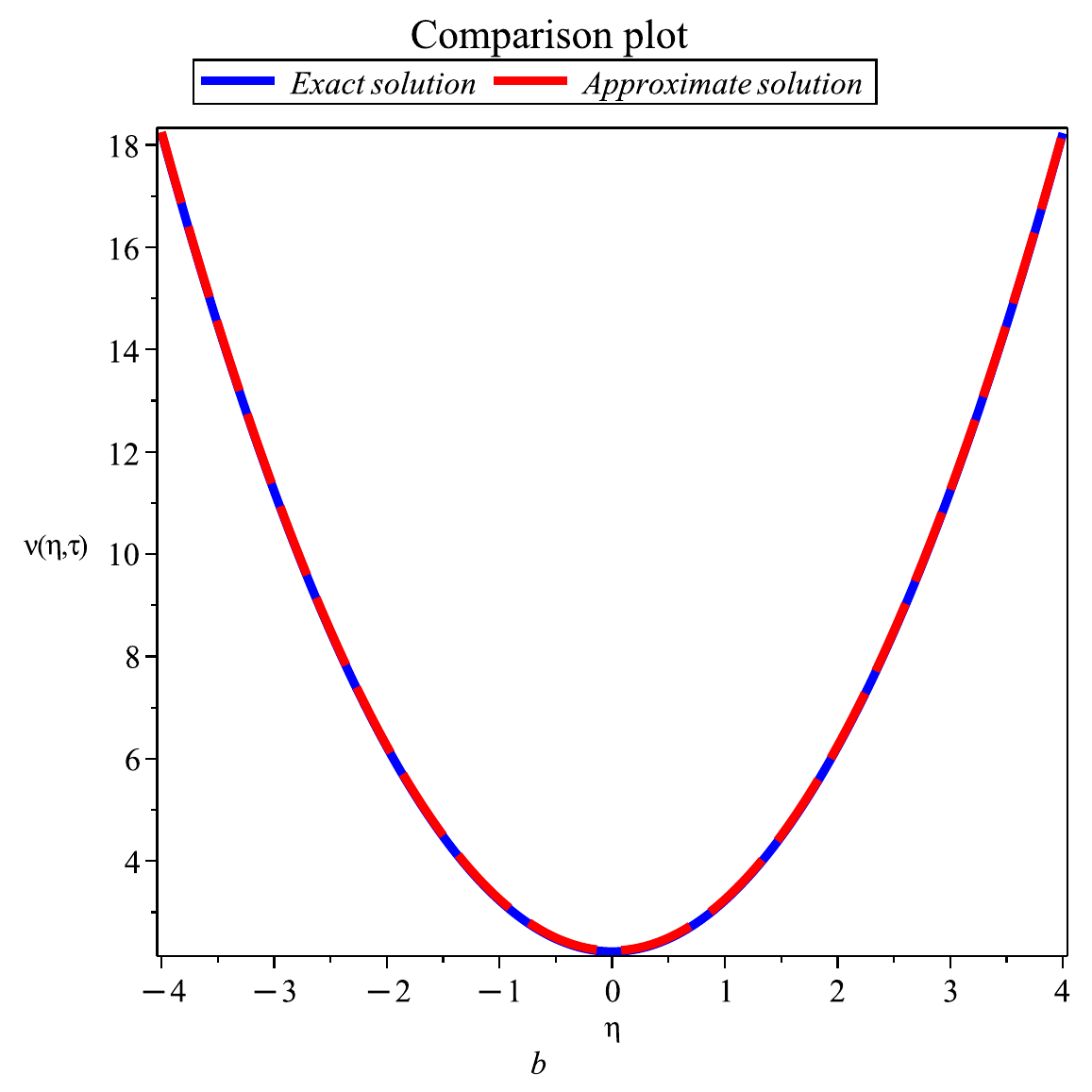}
	\caption{Comparisons  plots of approximate and exact solution for $(a)=u({{{x}}},{t})$ and $(b)=v({{{x}}},{t})$} of {example} (\ref{example32}).\label{figure3}
\end{figure}
\begin{figure}[H]\centering
	\includegraphics[width=8.0 cm]{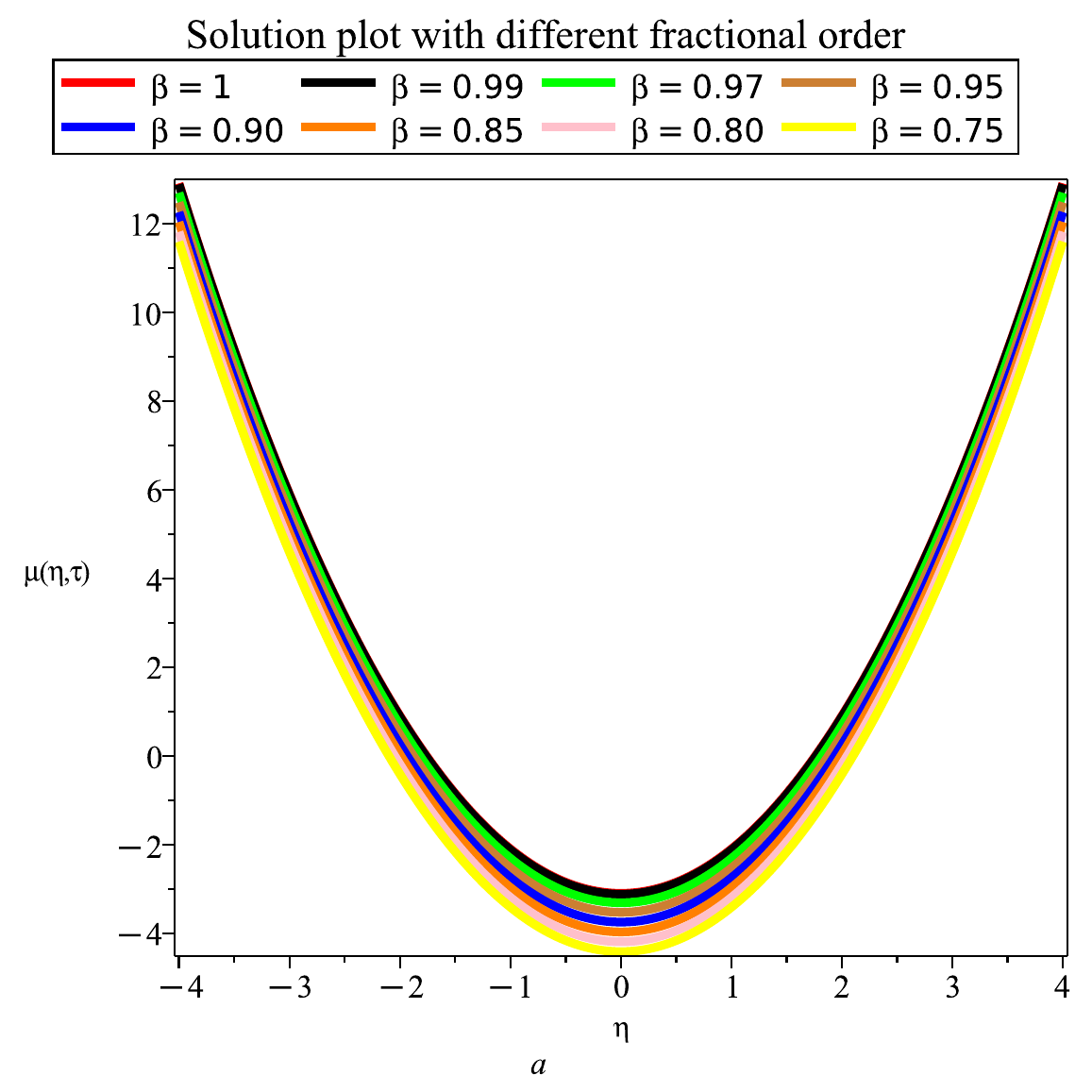}
	\includegraphics[width=8.0 cm]{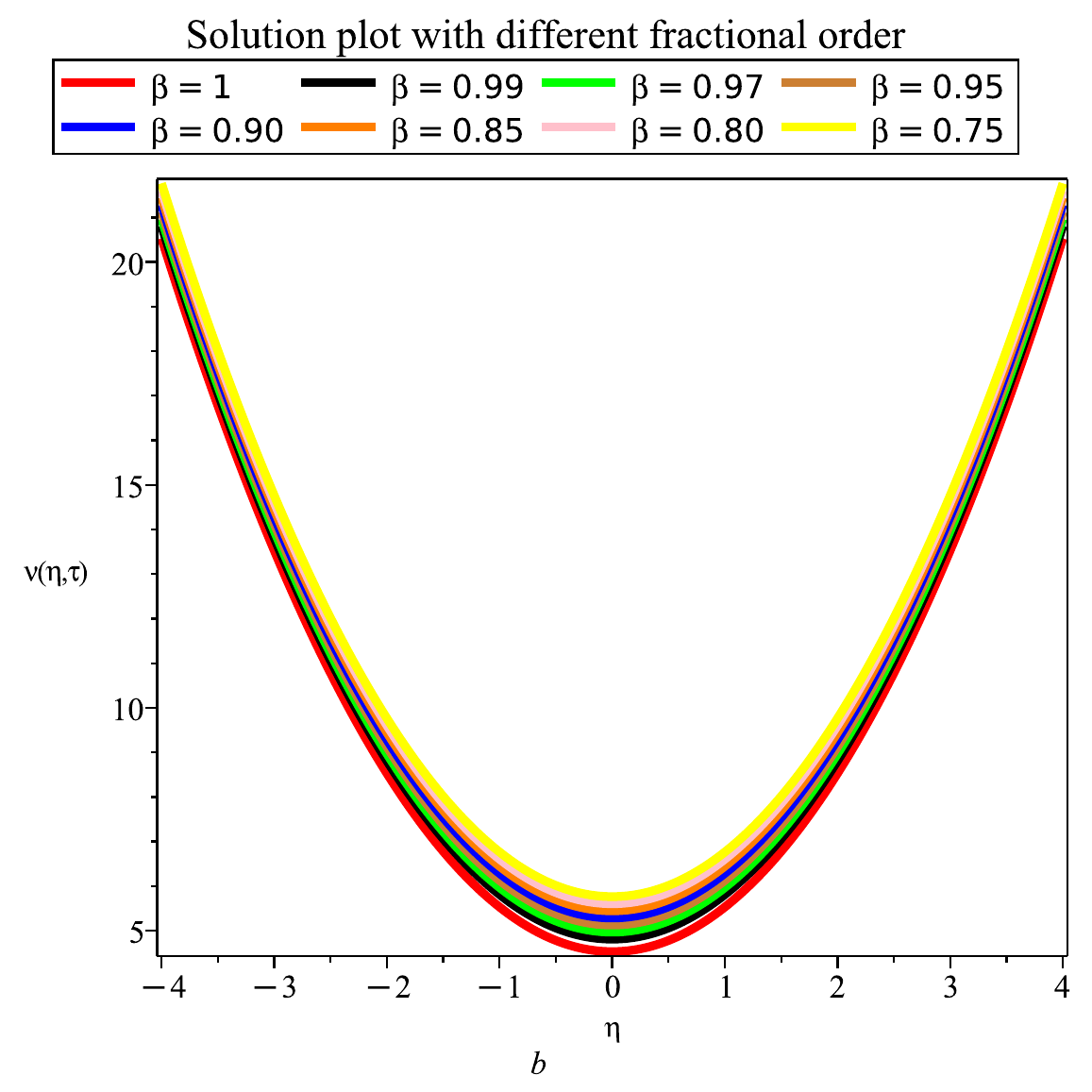}\\
\includegraphics[width=8.0 cm]{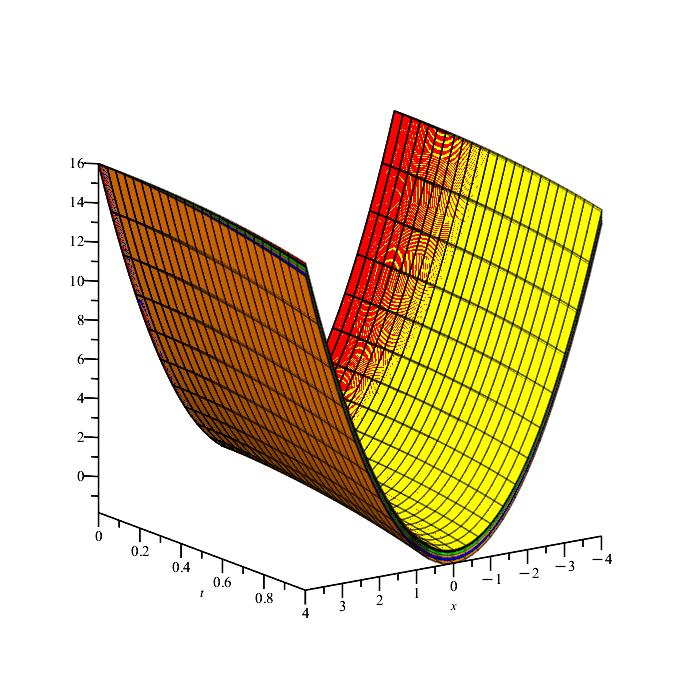}
\includegraphics[width=8.0 cm]{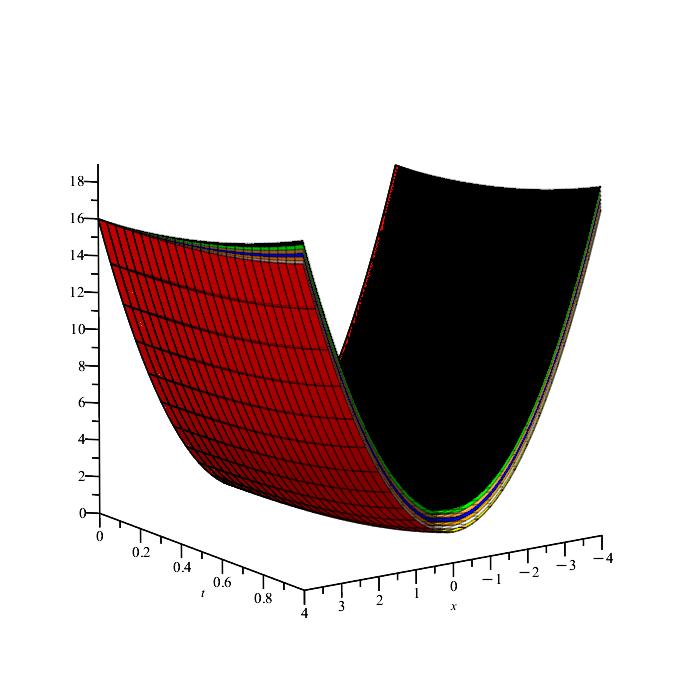}\\
\includegraphics[width=8.0 cm]{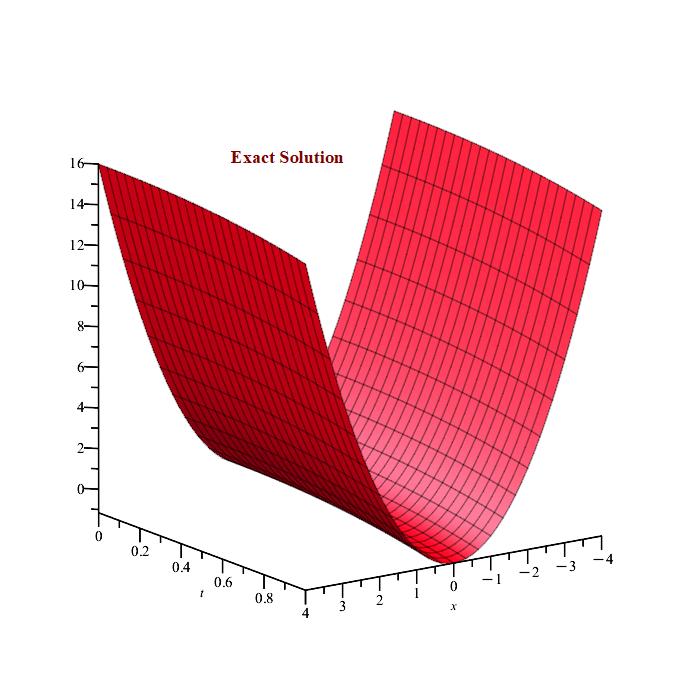}
\includegraphics[width=8.0 cm]{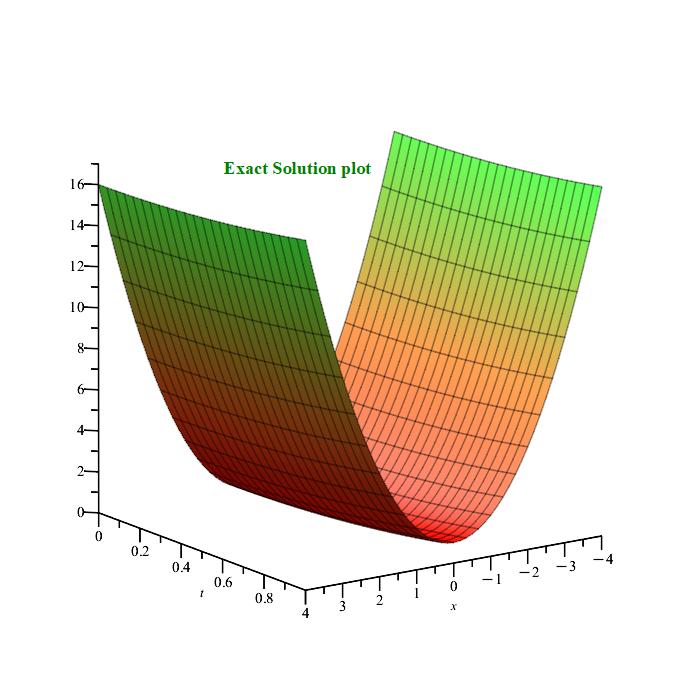}\\
	\caption{Comparisons 2D and 3D plotting for $u({{{x}}},{t})$ and $v({{{x}}},{t})$ with fractional order $\beta$ of  {example} (\ref{example32}).}\label{figure3}
\end{figure}

\begin{table}[H]
	\caption{\textbf{Numerical comparison of exact and approximate solutions with absolute error (AE) for various   time levels ${{t}}$  and different fractional  orders $\beta$ of example (\ref{example32}).}\label{tab4}}
	\centering
	\resizebox{\textwidth}{!}
	{
		\begin{tabular}{cccccccc}
			\toprule
			\textbf{{Variable}} & \textbf{time level} & \textbf{Exact } 	& \textbf{Approximate }	& \textbf{Absolute error}     & \textbf{Absolute error}& \textbf{Absolute error}&\textbf{Absolute error}\\
			\textbf{${u}$, ${v}$} & \textbf{${{t}}$} & \textbf{ solutions} 	& \textbf{solutions }	& \textbf{ $\beta=0.97$ }     & \textbf{ at $\beta=0.98$}& \textbf{ at $ \beta=0.99$}&\textbf{at $\beta=1$}\\
			\midrule  \midrule
 &0.01  &8.9999000000   &8.9999000000   &0.0000180170   &0.0000116810   &0.0000734600   &0.0000000000\\
 &0.02  &8.9996000000   &8.9996000000   &0.0000623480   &0.0000405720   &0.0002473550   &0.0000000000\\
 &0.03  &8.9991000000   &8.9991000000   &0.0001276930   &0.0000832730   &0.0004986010   &0.0000000000\\
 &0.04  &8.9983999990   &8.9984000000   &0.0002112820   &0.0001379930   &0.0008157370   &0.0000000010\\
 &0.05  &8.9974999990   &8.9975000000   &0.0003111910   &0.0002034890   &0.0011910130   &0.0000000010\\
 &0.06  &8.9964000000   &8.9964000000   &0.0004259440   &0.0002787980   &0.0016185820   &0.0000000000\\
 &0.07  &8.9951000000   &8.9951000000   &0.0005543400   &0.0003631390   &0.0020937450   &0.0000000000\\
 &0.08  &8.9935999990   &8.9936000000   &0.0006953610   &0.0004558470   &0.0026125790   &0.0000000010\\
 &0.09  &8.9919000000   &8.9919000000   &0.0008481270   &0.0005563490   &0.0031717170   &0.0000000000\\
 &0.10  &8.9899999990   &8.9900000000   &0.0010118580   &0.0006641350   &0.0037682070   &0.0000000010\\
 &0.11  &8.9879000000   &8.9879000000   &0.0011858600   &0.0007787450   &0.0043994210   &0.0000000000\\
 &0.12  &8.9856000000   &8.9856000000   &0.0013694930   &0.0008997680   &0.0050629860   &0.0000000000\\
 &0.13  &8.9831000000   &8.9831000000   &0.0015621810   &0.0010268190   &0.0057567480   &0.0000000000\\
 &0.14  &8.9803999990   &8.9804000000   &0.0017633820   &0.0011595470   &0.0064787200   &0.0000000010\\
 &0.15  &8.9775000000   &8.9775000000   &0.0019725990   &0.0012976230   &0.0072270720   &0.0000000000\\
 &0.16  &8.9744000000   &8.9744000000   &0.0021893650   &0.0014407410   &0.0080000950   &0.0000000000\\
\textbf{$u({{{x}}},{t})$} &0.17  &8.9711000000   &8.9711000000   &0.0024132410   &0.0015886120   &0.0087961970   &0.0000000000\\
 &0.18  &8.9676000000   &8.9676000000   &0.0026438090   &0.0017409600   &0.0096138740   &0.0000000000\\
 &0.19  &8.9639000010   &8.9639000000   &0.0028806800   &0.0018975280   &0.0104517170   &0.0000000010\\
 &0.20  &8.9600000020   &8.9600000000   &0.0031234780   &0.0020580730   &0.0113083880   &0.0000000020\\
 &0.21  &8.9559000020   &8.9559000000   &0.0033718490   &0.0022223550   &0.0121826160   &0.0000000020\\
 &0.22  &8.9516000040   &8.9516000000   &0.0036254520   &0.0023901520   &0.0130731990   &0.0000000040\\
 &0.23  &8.9471000060   &8.9471000000   &0.0038839610   &0.0025612480   &0.0139789820   &0.0000000060\\
 &0.24  &8.9424000080   &8.9424000000   &0.0041470660   &0.0027354380   &0.0148988670   &0.0000000080\\
 &0.25  &8.9375000120   &8.9375000000   &0.0044144620   &0.0029125220   &0.0158318030   &0.0000000120\\
 &0.26  &8.9324000170   &8.9324000000   &0.0046858630   &0.0030923090   &0.0167767800   &0.0000000170\\
 &0.27  &8.9271000220   &8.9271000000   &0.0049609910   &0.0032746130   &0.0177328300   &0.0000000220\\
 &0.28  &8.9216000290   &8.9216000000   &0.0052395720   &0.0034592580   &0.0186990210   &0.0000000290\\
 &0.29  &8.9159000390   &8.9159000000   &0.0055213490   &0.0036460700   &0.0196744550   &0.0000000390\\
 &0.30  &8.9100000520   &8.9100000000   &0.0058060700   &0.0038348790   &0.0206582710   &0.0000000520\\
 &0.31  &8.9039000680   &8.9039000000   &0.0060934890   &0.0040255270   &0.0216496330   &0.0000000680\\
 &0.32  &8.8976000870   &8.8976000000   &0.0063833700   &0.0042178540   &0.0226477350   &0.0000000870\\
 &0.33  &8.8911001120   &8.8911000000   &0.0066754850   &0.0044117070   &0.0236518040   &0.0000001120\\
 &0.34  &8.8844001410   &8.8844000000   &0.0069696080   &0.0046069400   &0.0246610850   &0.0000001410\\
 &0.35  &8.8775001790   &8.8775000000   &0.0072655210   &0.0048034050   &0.0256748520   &0.0000001790\\
 \bottomrule
 &0.01  &9.0001000000   &9.0001000000   &0.0000178650   &0.0000115880   &0.0000725190   &0.0000000000\\
 &0.02  &9.0003999990   &9.0004000000   &0.0000613120   &0.0000399320   &0.0002414110   &0.0000000010\\
 &0.03  &9.0009000000   &9.0009000000   &0.0001245400   &0.0000813130   &0.0004812930   &0.0000000000\\
 &0.04  &9.0016000010   &9.0016000000   &0.0002043670   &0.0001336800   &0.0007789800   &0.0000000010\\
 &0.05  &9.0025000000   &9.0025000000   &0.0002985250   &0.0001955600   &0.0011253210   &0.0000000000\\
 &0.06  &9.0036000000   &9.0036000000   &0.0004052290   &0.0002657930   &0.0015132780   &0.0000000000\\
 &0.07  &9.0049000000   &9.0049000000   &0.0005230040   &0.0003434180   &0.0019371320   &0.0000000000\\
 &0.08  &9.0063999990   &9.0064000000   &0.0006505820   &0.0004276110   &0.0023920680   &0.0000000010\\
 &0.09  &9.0080999990   &9.0081000000   &0.0007868640   &0.0005176450   &0.0028739300   &0.0000000010\\
 &0.10  &9.0100000010   &9.0100000000   &0.0009308630   &0.0006128810   &0.0033790760   &0.0000000010\\
 &0.11  &9.0121000000   &9.0121000000   &0.0010817000   &0.0007127350   &0.0039042590   &0.0000000000\\
 &0.12  &9.0144000000   &9.0144000000   &0.0012385650   &0.0008166810   &0.0044465650   &0.0000000000\\
 &0.13  &9.0169000010   &9.0169000000   &0.0014007270   &0.0009242320   &0.0050033470   &0.0000000010\\
 &0.14  &9.0196000000   &9.0196000000   &0.0015675030   &0.0010349400   &0.0055721930   &0.0000000000\\
 &0.15  &9.0225000010   &9.0225000000   &0.0017382660   &0.0011483950   &0.0061508850   &0.0000000010\\
 &0.16  &9.0256000000   &9.0256000000   &0.0019124290   &0.0012642010   &0.0067373790   &0.0000000000\\
 \textbf{$v({{{x}}},{t})$}&0.17  &9.0288999990   &9.0289000000   &0.0020894420   &0.0013820020   &0.0073297850   &0.0000000010\\
 &0.18  &9.0323999990   &9.0324000000   &0.0022687940   &0.0015014540   &0.0079263410   &0.0000000010\\
 &0.19  &9.0360999990   &9.0361000000   &0.0024500040   &0.0016222390   &0.0085254110   &0.0000000010\\
 &0.20  &9.0399999980   &9.0400000000   &0.0026326140   &0.0017440540   &0.0091254600   &0.0000000020\\
 &0.21  &9.0440999970   &9.0441000000   &0.0028161970   &0.0018666140   &0.0097250600   &0.0000000030\\
 &0.22  &9.0483999960   &9.0484000000   &0.0030003450   &0.0019896500   &0.0103228640   &0.0000000040\\
 &0.23  &9.0528999940   &9.0529000000   &0.0031846750   &0.0021129040   &0.0109176090   &0.0000000060\\
 &0.24  &9.0575999920   &9.0576000000   &0.0033688210   &0.0022361320   &0.0115081080   &0.0000000080\\
 &0.25  &9.0624999890   &9.0625000000   &0.0035524370   &0.0023591040   &0.0120932430   &0.0000000110\\
 &0.26  &9.0675999830   &9.0676000000   &0.0037351960   &0.0024816000   &0.0126719600   &0.0000000170\\
 &0.27  &9.0728999780   &9.0729000000   &0.0039167780   &0.0026034090   &0.0132432620   &0.0000000220\\
 &0.28  &9.0783999700   &9.0784000000   &0.0040968930   &0.0027243320   &0.0138062100   &0.0000000300\\
 &0.29  &9.0840999600   &9.0841000000   &0.0042752520   &0.0028441790   &0.0143599150   &0.0000000400\\
 &0.30  &9.0899999480   &9.0900000000   &0.0044515850   &0.0029627670   &0.0149035410   &0.0000000520\\
 &0.31  &9.0960999320   &9.0961000000   &0.0046256320   &0.0030799210   &0.0154362910   &0.0000000680\\
 &0.32  &9.1023999120   &9.1024000000   &0.0047971480   &0.0031954760   &0.0159574150   &0.0000000880\\
 &0.33  &9.1088998880   &9.1089000000   &0.0049658950   &0.0033092710   &0.0164662030   &0.0000001120\\
 &0.34  &9.1155998590   &9.1156000000   &0.0051316500   &0.0034211570   &0.0169619800   &0.0000001410\\
 &0.35  &9.1224998210   &9.1225000000   &0.0052941980   &0.0035309860   &0.0174441110   &0.0000001790\\

 \bottomrule
	\end{tabular}}
\end{table}

\subsection{Example}\label{example5}
Consider the non-linear coupled system of the thermo-elasticity
\begin{equation}\label{1dw4}
\begin{cases}
  \begin{split}
    &D_{{{{t}}}}^{{\beta+1}}{{{u}}}({{{{x}}}},{{{t}}})-\frac{1}{{{x}}}\left({{{x}}}{{{u}}}({{{{x}}}},{{{t}}})\frac{\partial}{\partial{{{{x}}}}}{{{u}}}({{{{x}}}},{{{t}}})\right)_{{{x}}}+{{{x}}}\frac{\partial {{{v}}}({{{{x}}}},{{{t}}})}{\partial{{{{x}}}}}+6{{{{x}}}}^{2}{{{t}}}^{2}=0,\ \  \ \ 0<\beta\leq1,\\
    &D_{{{{t}}}}^{{\beta}}{{{v}}}({{{{x}}}},{{{t}}})-\frac{1}{{{x}}}\left({{{x}}}{{{v}}}({{{{x}}}},{{{t}}})\frac{\partial}{\partial{{{{x}}}}}{{{v}}}({{{{x}}}},{{{t}}})\right)_{{{x}}}+{{{x}}}\frac{\partial^{2} {{{u}}}({{{{x}}}},{{{t}}})}{\partial{{{t}}}\partial{{{{x}}}}}-2{{{{x}}}}^{2}-2{{{{x}}}}^{2}{{{t}}}+8{{{{x}}}}^{2}{{{t}}}^4=0, \ \ 0<\beta\leq1,  \ \ {{{t}}}>0.
  \end{split}
  \end{cases}
\end{equation}

With initial conditions
\begin{equation}\label{1w5}
  {{{u}}}({{{{x}}}},0)=0,\ \ {{{u}}}_{{{{t}}}}({{{{x}}}},{{{t}}})={{{{x}}}}^{2}, \ \ {{{v}}}({{{{x}}}},0)=0.
\end{equation}
In system (\ref{1dw4}), we have two non-linear terms ${{{u}}}({{{{x}}}},{{{t}}})\frac{\partial}{\partial{{{{x}}}}}{{{u}}}({{{{x}}}},{{{t}}})$ and ${{{n}}}({{{{x}}}},{{{t}}})\frac{\partial}{\partial{{{{x}}}}}{{{v}}}({{{{x}}}},{{{t}}})$. Let ${\mathbf{N}_3({u},{v})}$ and $\mathbf{N}_4({u},{v})$ be the two non-linear terms. Now recall the definition \ref{ap} of the Adomian polynomial to control the non linearity i.e  ${\mathbf{N}_3({u},{v})}$ and ${\mathbf{N}_4({u},{v})}$.\\
Where
\begin{equation}
{\mathbf{N}_3({u},{v})}=\sum_{{{j}=0}}^{\infty}C_{{j}},\ \ \
{\mathbf{N}_4({u},{v})}=\sum_{{{j}=0}}^{\infty}D_{{j}}.
\end{equation}
\begin{equation*}
  \begin{split}
&C_{0}\coloneqq u_{0}\! \left({{{x}}} ,{t} \right) \left(\frac{\partial}{\partial x}u_{0}\! \left({{{x}}} ,{t} \right)\right)
 , \ \ D_{0}\coloneqq v_{0}\! \left({{{x}}} ,{t} \right) \left(\frac{\partial}{\partial x}v_{0}\! \left({{{x}}} ,{t} \right)\right)
,\\&C_{1}\coloneqq u_{1}\! \left({{{x}}} ,{t} \right) \left(\frac{\partial}{\partial x}u_{0}\! \left({{{x}}} ,{t} \right)\right)+u_{0}\! \left({{{x}}} ,{t} \right) \left(\frac{\partial}{\partial x}u_{1}\! \left({{{x}}} ,{t} \right)\right)
,\\&D_{1}\coloneqq v_{1}\! \left({{{x}}} ,{t} \right) \left(\frac{\partial}{\partial x}v_{0}\! \left({{{x}}} ,{t} \right)\right)+v_{0}\! \left({{{x}}} ,{t} \right) \left(\frac{\partial}{\partial x}v_{1}\! \left({{{x}}} ,{t} \right)\right)
,\\&C_{2}\coloneqq u_{2}\! \left({{{x}}} ,{t} \right) \left(\frac{\partial}{\partial x}u_{0}\! \left({{{x}}} ,{t} \right)\right)+u_{1}\! \left({{{x}}} ,{t} \right) \left(\frac{\partial}{\partial x}u_{1}\! \left({{{x}}} ,{t} \right)\right)+u_{0}\! \left({{{x}}} ,{t} \right) \left(\frac{\partial}{\partial x}u_{2}\! \left({{{x}}} ,{t} \right)\right)
,\\&D_{2}\coloneqq v_{2}\! \left({{{x}}} ,{t} \right) \left(\frac{\partial}{\partial x}v_{0}\! \left({{{x}}} ,{t} \right)\right)+v_{1}\! \left({{{x}}} ,{t} \right) \left(\frac{\partial}{\partial x}v_{1}\! \left({{{x}}} ,{t} \right)\right)+v_{0}\! \left({{{x}}} ,{t} \right) \left(\frac{\partial}{\partial x}v_{2}\! \left({{{x}}} ,{t} \right)\right)
,\\
\vdots
\\&C_{n}\coloneqq u_{n}\! \left({{{x}}} ,{t} \right) \left(\frac{\partial}{\partial x}u_{0}\! \left({{{x}}} ,{t} \right)\right)+u_{n-1}\! \left({{{x}}} ,{t} \right) \left(\frac{\partial}{\partial x}u_{1}\! \left({{{x}}} ,{t} \right)\right)+\cdots+u_{0}\! \left({{{x}}} ,{t} \right) \left(\frac{\partial}{\partial x}u_{n}\! \left({{{x}}} ,{t} \right)\right)
,\\&D_{n}\coloneqq v_{n}\! \left({{{x}}} ,{t} \right) \left(\frac{\partial}{\partial x}v_{0}\! \left({{{x}}} ,{t} \right)\right)+v_{n-1}\! \left({{{x}}} ,{t} \right) \left(\frac{\partial}{\partial x}v_{1}\! \left({{{x}}} ,{t} \right)\right)+\cdots+v_{0}\! \left({{{x}}} ,{t} \right) \left(\frac{\partial}{\partial x}v_{n}\! \left({{{x}}} ,{t} \right)\right), \forall \ \ {n} \ \ \in \ \ \mathbb{N}^{+}.
  \end{split}
\end{equation*}
Now, using the recurrence relation from the methodology section  \ref{s3}, particularly equation (\ref{129}) and equation (\ref{88}), we  obtained the following
\begin{equation*}
\begin{split}
&u_{0}\coloneqq {{{x}}}^{2} t \\
&v_{0}\coloneqq 0,\\
&u_{1}({{{x}}},{{t}})\coloneqq \frac{4 {{{x}}}^{2} {{t}}^{3+\beta}}{\Gamma \! \left(4+\beta \right)}
,\\&v_{1}({{{x}}},{{t}})\coloneqq \frac{2 {{{x}}}^{2} {{t}}^{\beta +1} \left(\beta^{3}-96 {{t}}^{3}+9 \beta^{2}+26 \beta +24\right)}{\Gamma \! \left(5+\beta \right)}
z\\&u_{2}({{{x}}},{{t}})\coloneqq 4 {{t}}^{2+2 \beta} {{{x}}}^{2} \left(-\frac{1}{\Gamma \! \left(3+2 \beta \right)}+\frac{16 {{t}}^{3} \left(10+\beta \right)}{\Gamma \! \left(6+2 \beta \right)}\right)
,\\&v_{2}({{{x}}},{{t}})\coloneqq -\frac{8 {{{x}}}^{2} {{t}}^{2+2 \beta}}{\Gamma \! \left(3+2 \beta \right)},\\
&u_{3}({{{x}}},{{t}})\coloneqq 16 {{t}}^{3+3 \beta} {{{x}}}^{2} \left(\frac{1}{\Gamma \! \left(4+3 \beta \right)}-\frac{4 {t} \left(3+2 \beta \right)}{\Gamma \! \left(5+3 \beta \right)}+\frac{8 {{t}}^{4} \left(\Gamma \! \left(7+2 \beta \right)+16 \Gamma \! \left(4+\beta \right)^{2} \left(10+\beta \right) \left(3+\beta \right)\right)}{\Gamma \! \left(4+\beta \right)^{2} \Gamma \! \left(8+3 \beta \right)}\right)
,\\&v_{3}({{{x}}},{{t}})\coloneqq 8 {{t}}^{1+3 \beta} {{{x}}}^{2} \left(\frac{\left(2+3 \beta \right) \Gamma \! \left(5+\beta \right)^{2}+4 t \Gamma \! \left(3+2 \beta \right) \left(2+\beta \right)^{2} \left(3+\beta \right)^{2} \left(4+\beta \right)^{2}}{\Gamma \! \left(5+\beta \right)^{2} \Gamma \! \left(3+3 \beta \right)}-\frac{16 {{t}}^{3} \left(10+\beta \right)}{\Gamma \! \left(5+3 \beta \right)}\right),\\
\vdots
\end{split}
\end{equation*}
It is clear that the decomposition series form solution for ${{u}}({{{{x}}}},{{{t}}})$ and ${{v}}({{{{x}}}},{{{t}}})$ is given by,
\begin{equation}\label{1w4}
\begin{split}
 & {u}({{{x}}},{{t}})=\sum _{n=0}^{\infty }{u}_{{n}} ({{{x}}},{{t}})={u}_{0}({{{x}}},{{t}})+{u}_{1}({{{x}}},{{t}})+{u}_{2}({{{x}}},{{t}})+{u}_{3}({{{x}}},{{t}})+\ldots\\
 &{v}({{{x}}},{{t}})=\sum _{n=0}^{\infty }{v}_{{n}} ({{{x}}},{{t}})={v}_{0}({{{x}}},{{t}})+{v}_{1}({{{x}}},{{t}})+{v}_{2}({{{x}}},{{t}})+{v}_{3}({{{x}}},{{t}})+\ldots\\
 \end{split}
\end{equation}
Putting $u_{0}, \  u_{1},\ \cdots,$ and $v_{0}, \ v_{1}, \ \cdots,$ in equation (\,\ref{1w4}\,) and setting $\beta=1$, we obtain

\begin{equation}\label{1q5}
\begin{split}
u({{{x}}},{t})=&{{{x}}}^{2}{{t}} +\frac{{{{x}}}^{2} {{t}}^{4}}{6}+\frac{{{{x}}}^{2} {{t}}^{4} \left(88 {{t}}^{3}-105\right)}{630}+\frac{{{{x}}}^{2} {{t}}^{6} \left(43 {{t}}^{4}-100{{t}} +35\right)}{1575}\\
&-\frac{{{{x}}}^{2} {{t}}^{6} \left(37868 {{t}}^{7}-111254 {{t}}^{4}-25740 {{t}}^{3}+154440{{t}} +45045\right)}{2027025}
+\ldots~~~~ \longrightarrow {{{{x}}}}^{2}{{{t}}},\\
v({{{x}}},{t})=&  {{{x}}}^{2} \left(-\frac{8}{5} {{t}}^{5}+{{t}}^{2}\right)-\frac{{{{x}}}^{2} {{t}}^{4}}{3}+\frac{{{{x}}}^{2} {{t}}^{4} \left(32256 {{t}}^{7}-55440 {{t}}^{4}-4840 {{t}}^{3}+27720{{t}} +5775\right)}{17325}\\&+\frac{2 {{{x}}}^{2} {{t}}^{6} \left(629 {{t}}^{4}-500 t -{t}35\right)}{1575}+\ldots~~~  \longrightarrow {{{{x}}}}^{2}{{{t}}}^{2}.
\end{split}
\end{equation}
Implies that
$ {{{u}}}({{{{x}}}},{{{t}}})={{{{x}}}}^{2}{{{t}}}$ and ${{{v}}}({{{{x}}}},{{{t}}})={{{{x}}}}^{2}{{{t}}}^{2}$. Which  is an exact solution of  for problem (\ref{1dw4}). In addition, we observe that the final result can be compared to the outcome achieved by LADM (see \cite{thermo32a}).
\begin{equation}\label{16}
  \begin{split}
     {{{u}}}({{_{{{x}}}}},{{{t}}})={{{{x}}}}^{2}{{{t}}},\\
     {{{v}}}({{{{x}}}},{{{t}}})={{{{x}}}}^{2}{{{t}}}^{2}.
  \end{split}
\end{equation}
\begin{figure}[H]\centering
	\includegraphics[width=8.0 cm]{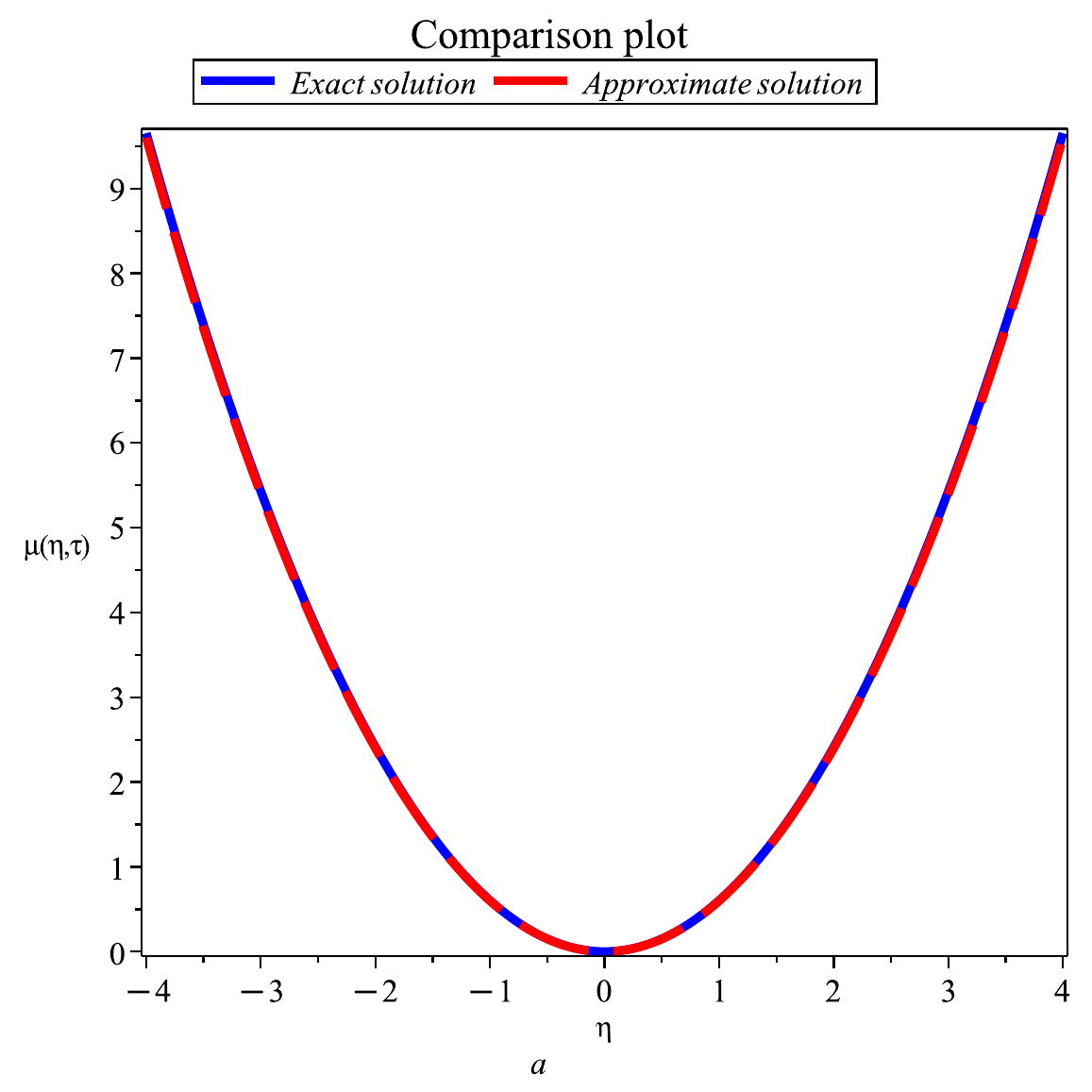}
	\includegraphics[width=8.0 cm]{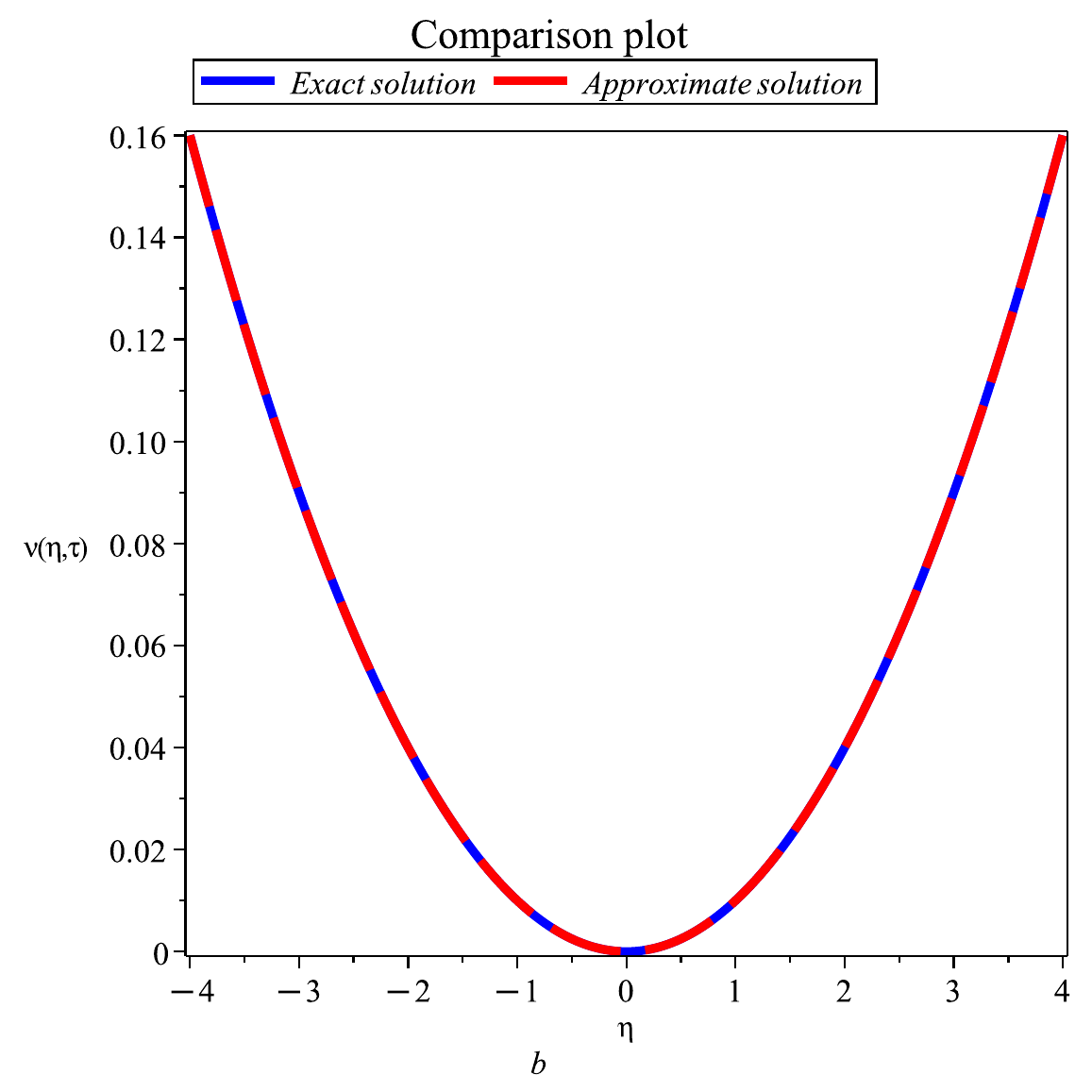}
	\caption{Comparisons  plots of approximate and exact solution for $(a)=u({{{x}}},{t})$ and $(b)=v({{{x}}},{t})$} of {example} (\ref{example5}).\label{figure3}
\end{figure}
\begin{figure}[H]\centering
	\includegraphics[width=8.0 cm]{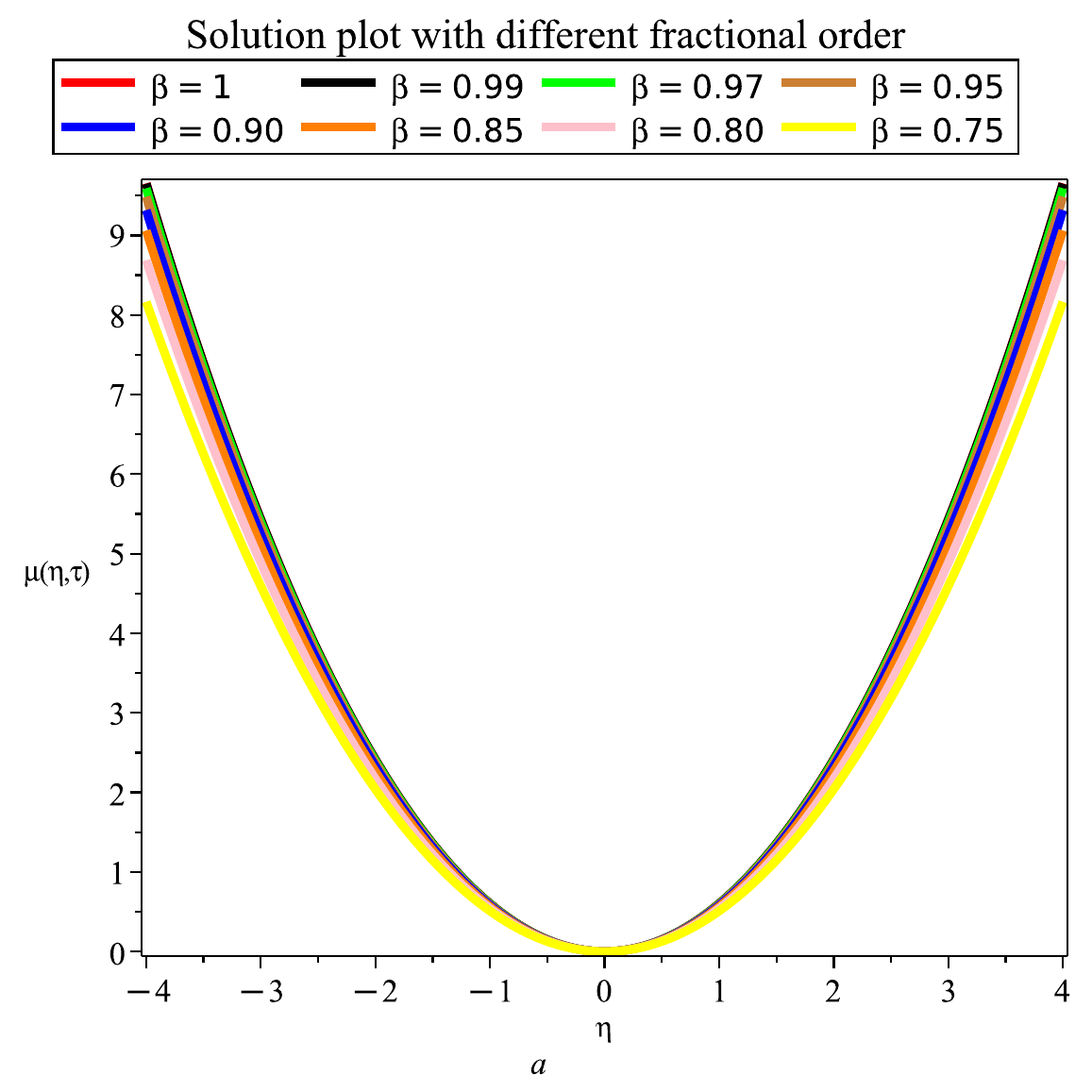}
	\includegraphics[width=8.0 cm]{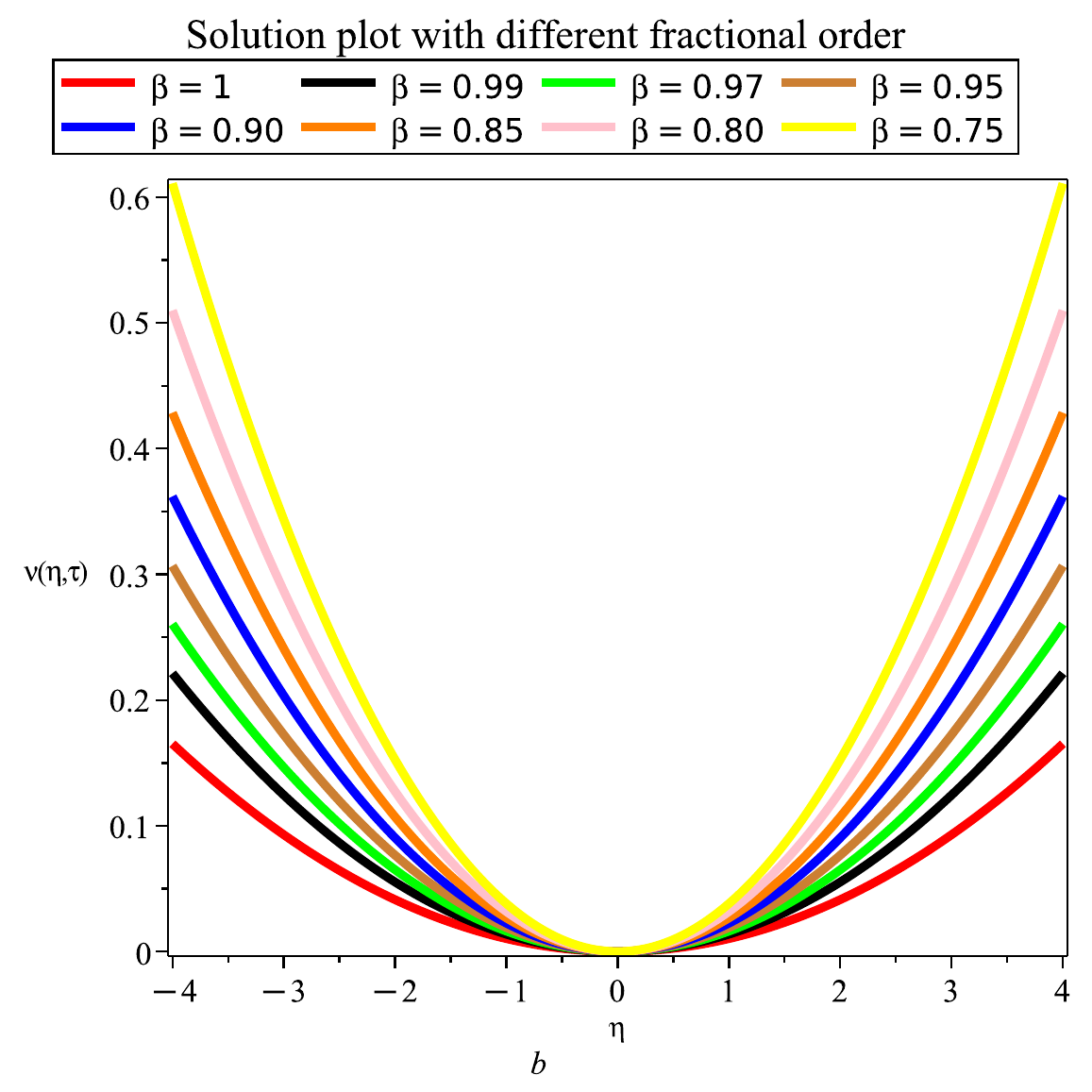}\\
\includegraphics[width=8.0 cm]{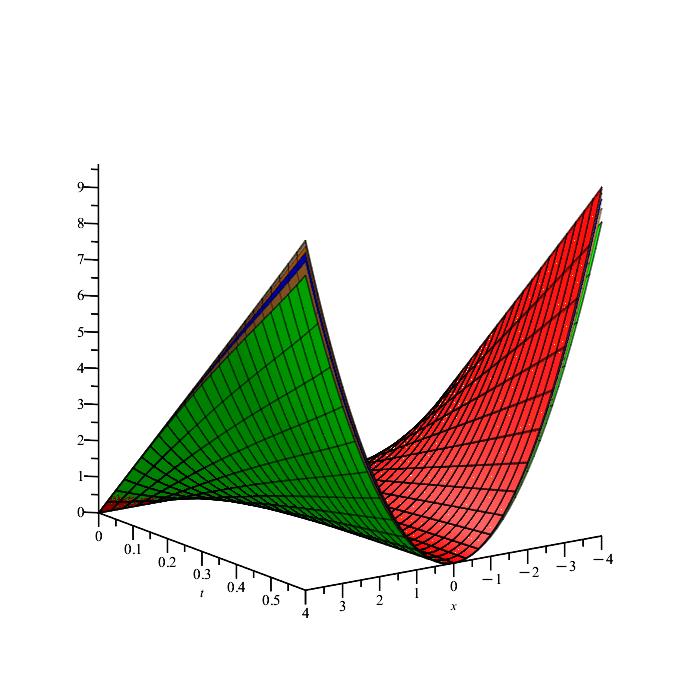}
\includegraphics[width=8.0 cm]{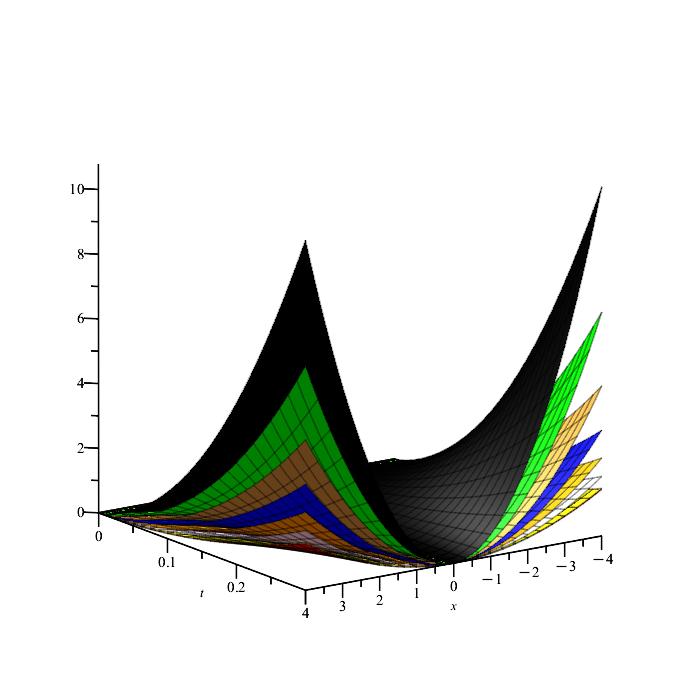}\\
\includegraphics[width=8.0 cm]{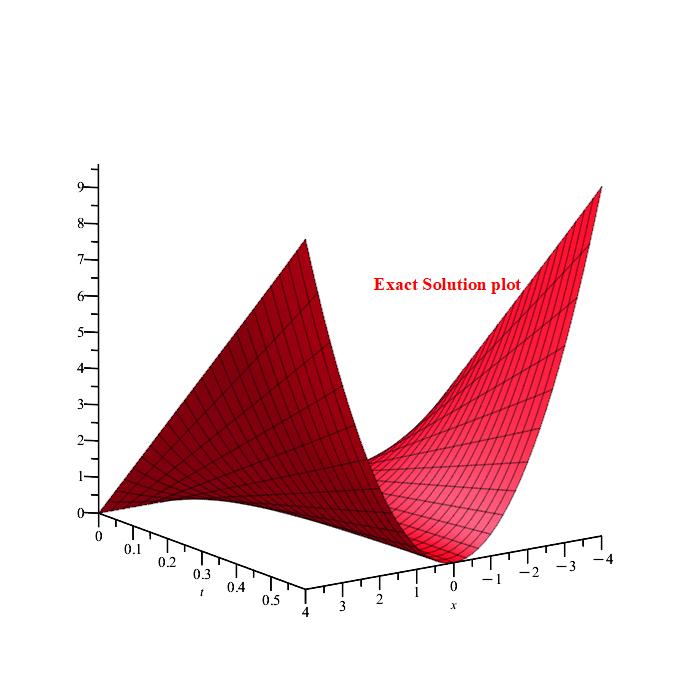}
\includegraphics[width=8.0 cm]{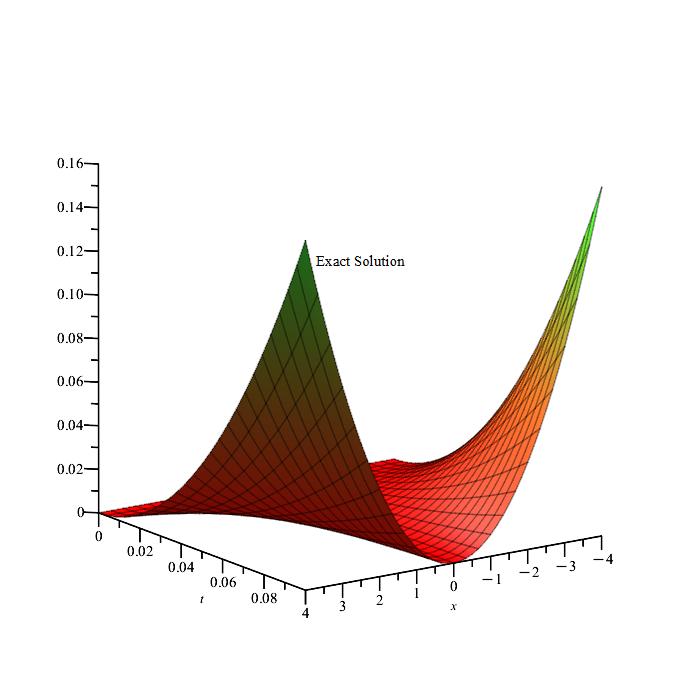}\\
	\caption{Comparisons 2D and 3D plotting for $u({{{x}}},{t})$ and $v({{{x}}},{t})$ with fractional order $\beta$ of  {example} (\ref{example5}).}\label{figure3}
\end{figure}

\begin{table}[H]
	\caption{\textbf{Numerical comparison of exact and approximate solutions with absolute error at different time levels ${{t}}$  and fractional  orders $\beta$ of example (\ref{example5}).} \label{tab5}}
	\centering
	\resizebox{\textwidth}{!}
	{
		\begin{tabular}{cccccccc}
			\toprule
			\textbf{{Variable}} & \textbf{time level} & \textbf{Exact } 	& \textbf{Approximate }	& \textbf{Absolute error}     & \textbf{Absolute error}& \textbf{Absolute error}&\textbf{Absolute error}\\
			\textbf{${u}$, ${v}$} & \textbf{${{t}}$} & \textbf{ solutions} 	& \textbf{solutions }	& \textbf{ $\beta=0.97$ }     & \textbf{ at $\beta=0.98$}& \textbf{ at $ \beta=0.99$}&\textbf{at $\beta=1$}\\
			\midrule  \midrule
 &0.01  &0.0900000000   &0.0900000000   &0.0000000036   &0.0000000022   &0.0000000231   &0.0000000000\\
 &0.02  &0.1800000000   &0.1800000000   &0.0000000498   &0.0000000305   &0.0000002939   &0.0000000000\\
 &0.03  &0.2700000001   &0.2700000000   &0.0000002284   &0.0000001412   &0.0000012920   &0.0000000001\\
 &0.04  &0.3600000009   &0.3600000000   &0.0000006709   &0.0000004164   &0.0000036789   &0.0000000009\\
 &0.05  &0.4500000037   &0.4500000000   &0.0000015423   &0.0000009602   &0.0000082605   &0.0000000037\\
 &0.06  &0.5400000112   &0.5400000000   &0.0000030375   &0.0000018947   &0.0000159631   &0.0000000112\\
 &0.07  &0.6300000292   &0.6300000000   &0.0000053754   &0.0000033577   &0.0000278150   &0.0000000292\\
 &0.08  &0.7200000668   &0.7200000000   &0.0000087964   &0.0000054988   &0.0000449316   &0.0000000668\\
 &0.09  &0.8100001391   &0.8100000000   &0.0000135567   &0.0000084768   &0.0000685029   &0.0000001391\\
 &0.10  &0.9000002686   &0.9000000000   &0.0000199256   &0.0000124558   &0.0000997802   &0.0000002686\\
 &0.11  &0.9900004880   &0.9900000000   &0.0000281797   &0.0000176020   &0.0001400659   &0.0000004880\\
 &0.12  &1.0800008430   &1.0800000000   &0.0000385990   &0.0000240790   &0.0001907020   &0.0000008430\\
 &0.13  &1.1700013960   &1.1700000000   &0.0000514610   &0.0000320440   &0.0002530570   &0.0000013960\\
 &0.14  &1.2600022290   &1.2600000000   &0.0000670340   &0.0000416390   &0.0003285190   &0.0000022290\\
 &0.15  &1.3500034510   &1.3500000000   &0.0000855760   &0.0000529950   &0.0004184790   &0.0000034510\\
 &0.16  &1.4400051990   &1.4400000000   &0.0001073180   &0.0000662120   &0.0005243210   &0.0000051990\\
 &0.17  &1.5300076470   &1.5300000000   &0.0001324690   &0.0000813630   &0.0006474130   &0.0000076470\\
\textbf{$u({{{x}}},{t})$} &0.18  &1.6200110090   &1.6200000000   &0.0001611950   &0.0000984790   &0.0007890880   &0.0000110090\\
 &0.19  &1.7100155530   &1.7100000000   &0.0001936200   &0.0001175490   &0.0009506310   &0.0000155530\\
 &0.20  &1.8000216020   &1.8000000000   &0.0002298110   &0.0001385040   &0.0011332680   &0.0000216020\\
 &0.21  &1.8900295440   &1.8900000000   &0.0002697710   &0.0001612070   &0.0013381490   &0.0000295440\\
 &0.22  &1.9800398460   &1.9800000000   &0.0003134230   &0.0001854480   &0.0015663250   &0.0000398460\\
 &0.23 & 2.0700530570  & 2.0700000000   &0.0003606000   &0.0002109250   &0.0018187450   &0.0000530570\\
 &0.24  &2.1600698260  & 2.1600000000   &0.0004110340   &0.0002372360   &0.0020962190   &0.0000698260\\
 &0.25  &2.2500909150  & 2.2500000000   &0.0004643360   &0.0002638680   &0.0023994120   &0.0000909150\\
 &0.26 & 2.3401172050  & 2.3400000000   &0.0005199880   &0.0002901710   &0.0027288170   &0.0001172050\\
 &0.27  &2.4301497190   &2.4300000000   &0.0005773190   &0.0003153570   &0.0030847320   &0.0001497190\\
 &0.28 & 2.5201896280  & 2.5200000000   &0.0006354890   &0.0003384680   &0.0034672390   &0.0001896280\\
 &0.29  &2.6102382840  & 2.6100000000   &0.0006934720   &0.0003583680   &0.0038761730   &0.0002382840\\
 &0.30 & 2.7002972160  & 2.7000000000   &0.0007500330   &0.0003737200   &0.0043111030   &0.0002972160\\
 &0.31 & 2.7903681730  & 2.7900000000   &0.0008037090   &0.0003829630   &0.0047712900   &0.0003681730\\
 &0.32  &2.8804531240  & 2.8800000000   &0.0008527790   &0.0003842930   &0.0052556710   &0.0004531240\\
 &0.33  &2.9705542980  & 2.9700000000   &0.0008952460   &0.0003756310   &0.0057628120   &0.0005542980\\
 &0.34  &3.0606741940  & 3.0600000000   &0.0009288040   &0.0003546090   &0.0062908850   &0.0006741940\\
 &0.35 & 3.1508156150  & 3.1500000000   &0.0009508140   &0.0003185290   &0.0068376190   &0.0008156150\\
  \bottomrule
 &0.01  &0.0009000000   &0.0009000000   &0.0001621642   &0.0001051379   &0.0006612578   &0.0000000000\\
 &0.02  &0.0036000000   &0.0036000000   &0.0005613200   &0.0003652578   &0.0022276244   &0.0000000000\\
 &0.03  &0.0080999999   &0.0081000000   &0.0011501279   &0.0007499863   &0.0044937882   &0.0000000001\\
 &0.04  &0.0143999994   &0.0144000000   &0.0019042450   &0.0012435780   &0.0073602004   &0.0000000006\\
 &0.05  &0.0224999969   &0.0225000000   &0.0028072064   &0.0018353317   &0.0107618910   &0.0000000031\\
 &0.06  &0.0323999881   &0.0324000000   &0.0038467892   &0.0025172689   &0.0146521293   &0.0000000119\\
 &0.07  &0.0440999627   &0.0441000000   &0.0050134902   &0.0032831592   &0.0189957314   &0.0000000373\\
 &0.08  &0.0575998990   &0.0576000000   &0.0062997374   &0.0041280126   &0.0237656498   &0.0000001010\\
 &0.09  &0.0728997558   &0.0729000000   &0.0076994250   &0.0050477779   &0.0289409934   &0.0000002442\\
 &0.10  &0.0899994607   &0.0900000000   &0.0092076106   &0.0060391442   &0.0345057680   &0.0000005393\\
 &0.11  &0.1088988932   &0.1089000000   &0.0108203021   &0.0070993977   &0.0404480156   &0.0000011068\\
 &0.12  &0.1295978620   &0.1296000000   &0.0125342939   &0.0082263093   &0.0467591862   &0.0000021380\\
 &0.13  &0.1520960760   &0.1521000000   &0.0143470322   &0.0094180361   &0.0534336513   &0.0000039240\\
 &0.14  &0.1763931062   &0.1764000000   &0.0162564938   &0.0106730316   &0.0604683040   &0.0000068938\\
 &0.15  &0.2024883375   &0.2025000000   &0.0182610701   &0.0119899529   &0.0678622098   &0.0000116625\\
 &0.16  &0.2303809106   &0.2304000000   &0.0203594486   &0.0133675643   &0.0756162878   &0.0000190894\\
\textbf{$v({{{x}}},{t})$} &0.17  &0.2600696502   &0.2601000000   &0.0225504882   &0.0148046328   &0.0837329975   &0.0000303498\\
 &0.18  &0.2915529772   &0.2916000000   &0.0248330824   &0.0162998074   &0.0922160317   &0.0000470228\\
 &0.19  &0.3248288082   &0.3249000000   &0.0272060085   &0.0178514865   &0.1010699922   &0.0000711918\\
 &0.20  &0.3598944323   &0.3600000000   &0.0296677571   &0.0194576663   &0.1103000505   &0.0001055677\\
 &0.21  &0.3967463724   &0.3969000000   &0.0322163411   &0.0211157671   &0.1199115834   &0.0001536276\\
 &0.22  &0.4353802214   &0.4356000000   &0.0348490809   &0.0228224366   &0.1299097759   &0.0002197786\\
 &0.23  &0.4757904538   &0.4761000000   &0.0375623620   &0.0245733280   &0.1402991928   &0.0003095462\\
 &0.24  &0.5179702160   &0.5184000000   &0.0403513630   &0.0263628500   &0.1510833060   &0.0004297840\\
 &0.25  &0.5619110821   &0.5625000000   &0.0432097523   &0.0281838848   &0.1622639830   &0.0005889179\\
 &0.26  &0.6076027854   &0.6084000000   &0.0461293503   &0.0300274746   &0.1738409261   &0.0007972146\\
 &0.27  &0.6550329148   &0.6561000000   &0.0490997533   &0.0318824711   &0.1858110610   &0.0010670852\\
 &0.28  &0.7041865763   &0.7056000000   &0.0521079214   &0.0337351508   &0.1981678763   &0.0014134237\\
 &0.29  &0.7550460229   &0.7569000000   &0.0551377212   &0.0355687880   &0.2109007038   &0.0018539771\\
 &0.30  &0.8075902411   &0.8100000000   &0.0581694295   &0.0373631871   &0.2239939430   &0.0024097589\\
 &0.31  &0.8617945033   &0.8649000000   &0.0611791902   &0.0390941750   &0.2374262270   &0.0031054967\\
 &0.32  &0.9176298758   &0.9216000000   &0.0641384234   &0.0407330422   &0.2511695340   &0.0039701242\\
 &0.33  &0.9750626838   &0.9801000000   &0.0670131890   &0.0422459480   &0.2651882150   &0.0050373162\\
 &0.34  &1.0340539380   &1.0404000000   &0.0697635020   &0.0435932650   &0.2794379980   &0.0063460620\\
 &0.35  &1.0945587060   &1.1025000000   &0.0723425910   &0.0447288900   &0.2938648940   &0.0079412940\\
 \bottomrule
	\end{tabular}}
\end{table}
\section{Conclusion}\label{s6}
In this article, an efficient and accurate technique is utilized to solve nonlinear fractional systems of thermo-elastic models. The non-linearity is addressed using Adomian polynomials, while the Aboodh transformation is employed to simplify the fractional derivatives in each targeted problem. Ultimately, the Adomian decomposition method is implemented to achieve the complete solution for each problem.
Tables \ref{table1}, \ref{table2}, \ref{table3}, \ref{tab4}, and \ref{tab5}
 present the exact, approximate, and associated absolute errors at different fractional orders. It is observed from these tables that as the fractional orders approach integer values, the solutions at fractional orders converge to the integer-order solutions. The figures illustrate a close relationship between the exact, approximate, and fractional-order solutions. Additionally, it is noted that the current procedure is straightforward and requires fewer calculations.
In the future, this methodology can be extended to solve other nonlinear systems of fractional partial differential equations.
 \section*{Declaration}
\subsection*{Conflict of interest}
	The authors declare that they have no known competing financial interests or personal relationships that could have appeared to influence the work reported in this research paper.
\subsection*{Ethical approval}
The conducted research is not related to either human or animals use.
\subsection*{Data availability statement }
Data sharing is not applicable to this article as no datasets were generated or analysed during the current study.
\bibliographystyle{unsrt}
\bibliography{ref_for_thermo_paper}
 \end{document}